%% file: main.tex
\documentclass[a4paper,UKenglish,cleveref, autoref, thm-restate]{lipics-v2021}
\usepackage{pythonhighlight}

\title{Column Generation with Domain-Independent Dynamic Programming} %TODO Please add

\author{Ryo Kuroiwa}{
Principles of Informatics Division, National Institute of Informatics, Tokyo, Japan \and The Graduate University for Advanced Studies, SOKENDAI, Kanagawa, Japan \and \url{https://researchmap.jp/rkuroiwa?lang=en} }{kuroiwa@nii.ac.jp}{https://orcid.org/0000-0002-3753-1644}{JSPS KAKENHI grant number JP25K24378}
\author{Edward Lam}{
Department of Data Science and Artificial Intelligence, Monash University, Melbourne, Australia \and \url{https://ed-lam.com/} }{edward.lam@monash.edu}{https://orcid.org/0000-0002-4485-5014}{Australian Research Council, grant DE240100042}

\authorrunning{R. Kuroiwa and E. Lam} %TODO mandatory. First: Use abbreviated first/middle names. Second (only in severe cases): Use first author plus 'et al.'

\Copyright{Ryo Kuroiwa and Edward Lam} %TODO mandatory, please use full first names. LIPIcs license is "CC-BY";  http://creativecommons.org/licenses/by/3.0/

\ccsdesc[500]{Computing methodologies~Modeling methodologies}
\ccsdesc[500]{Theory of computation~Dynamic programming}
\ccsdesc[500]{Theory of computation~Mathematical optimization}

\keywords{Modelling \& Modelling Languages,Dynamic Programming,Operations Research \& Mathematical Optimisation} %TODO mandatory; please add comma-separated list of keywords

\category{} %optional, e.g. invited paper

\relatedversion{} %optional, e.g. full version hosted on arXiv, HAL, or other respository/website
\supplement{}
\supplementdetails[linktext={Code used for experiments}]{Software}{https://github.com/Kurorororo/cp2026-cg-didp-code}

\nolinenumbers %uncomment to disable line numbering

\EventEditors{Nicolas Beldiceanu}
\EventNoEds{1}
\EventLongTitle{32nd International Conference on Principles and Practice of Constraint Programming (CP 2026)}
\EventShortTitle{CP 2026}
\EventAcronym{CP}
\EventYear{2026}
\EventDate{July 20--23, 2026}
\EventLocation{Lisbon, Portugal}
\EventLogo{}
\SeriesVolume{379}
\ArticleNo{36}
\usepackage{algorithm}
\usepackage[noend]{algpseudocode}
\usepackage{subfiles}

\usepackage{microtype}

\usepackage{soul}

\usepackage{amsmath}
\usepackage{amsfonts}          % Math fonts
\usepackage{amssymb}           % Extra math symbols
\usepackage{mathtools}         % Commands like \phantom for math
\usepackage{empheq}            % Text to the side of equations
\allowdisplaybreaks[1]         % Allow equation groups to break across pages
\mathtoolsset{centercolon}     % Make colon with middle vertical alignment

\crefname{prob}{Problem}{Problems}
\crefformat{prob}{#2Problem~(#1)#3}
\Crefformat{prob}{#2Problem~(#1)#3}
\crefrangeformat{prob}{#3Problems~(#1)#4 to~#5(#2)#6}
\Crefrangeformat{prob}{#3Problems~(#1)#4 to~#5(#2)#6}
\crefmultiformat{prob}{#2Problems~(#1)#3}{ and~#2(#1)#3}{, #2(#1)#3}{ and~#2(#1)#3}
\Crefmultiformat{prob}{#2Problems~(#1)#3}{ and~#2(#1)#3}{, #2(#1)#3}{ and~#2(#1)#3}

\crefname{constr}{Constraint}{Constraints}
\crefformat{constr}{#2Constraint~(#1)#3}
\Crefformat{constr}{#2Constraint~(#1)#3}
\crefrangeformat{constr}{#3Constraints~(#1)#4 to~#5(#2)#6}
\Crefrangeformat{constr}{#3Constraints~(#1)#4 to~#5(#2)#6}
\crefmultiformat{constr}{#2Constraints~(#1)#3}{ and~#2(#1)#3}{, #2(#1)#3}{ and~#2(#1)#3}
\Crefmultiformat{constr}{#2Constraints~(#1)#3}{ and~#2(#1)#3}{, #2(#1)#3}{ and~#2(#1)#3}

\crefname{ineq}{Inequality}{Inequalities}
\crefformat{ineq}{#2Inequality~(#1)#3}
\Crefformat{ineq}{#2Inequality~(#1)#3}
\crefrangeformat{ineq}{#3Inequalities~(#1)#4 to~#5(#2)#6}
\Crefrangeformat{ineq}{#3Inequalities~(#1)#4 to~#5(#2)#6}
\crefmultiformat{ineq}{#2Inequalities~(#1)#3}{ and~#2(#1)#3}{, #2(#1)#3}{ and~#2(#1)#3}
\Crefmultiformat{ineq}{#2Inequalities~(#1)#3}{ and~#2(#1)#3}{, #2(#1)#3}{ and~#2(#1)#3}

\crefname{target}{Objective}{Objectives}
\crefformat{target}{#2Objective~(#1)#3}
\Crefformat{target}{#2Objective~(#1)#3}
\crefrangeformat{target}{#3Objectives~(#1)#4 to~#5(#2)#6}
\Crefrangeformat{target}{#3Objectives~(#1)#4 to~#5(#2)#6}
\crefmultiformat{target}{#2Objectives~(#1)#3}{ and~#2(#1)#3}{, #2(#1)#3}{ and~#2(#1)#3}
\Crefmultiformat{target}{#2Objectives~(#1)#3}{ and~#2(#1)#3}{, #2(#1)#3}{ and~#2(#1)#3}

\usepackage{pgfplots}           % Plots
\pgfplotsset{compat=newest}     % Set version
\usepgfplotslibrary{groupplots} % Subplots
\definecolor{lightgrey}{RGB}{218,218,218}
\definecolor{blue}{HTML}{1f77b4}
\definecolor{orange}{HTML}{ff7f0e}
\definecolor{green}{HTML}{2ca02c}
\definecolor{red}{HTML}{d62728}
\definecolor{purple}{HTML}{9467bd}
\definecolor{brown}{HTML}{8c564b}
\definecolor{pink}{HTML}{e377c2}
\definecolor{grey}{HTML}{7f7f7f}
\definecolor{olive}{HTML}{bcbd22}
\definecolor{teal}{HTML}{17becf}

\tikzset{
    thinline/.style = {
        line width=0.1mm
    },
    regularline/.style = {
        line width=0.25mm
    },
    thickline/.style = {
        line width=0.3mm
    },
    plotline solid/.style = {
        regularline,
        solid
    },
    plotline dashed/.style = {
        regularline,
        dash pattern=on 3pt off 2pt
    },
    plotline dotted/.style = {
        regularline,
        dash pattern=on 0.7pt off 1.6pt,
        line cap=round
    },
    plotline dashdot/.style = {
        regularline,
        dash pattern=on 3pt off 1.5pt on 0.7pt off 1.5pt,
        line cap=round
    },
    plotline longdash/.style = {
        regularline,
        dash pattern=on 5pt off 2pt
    },
    plotline dashdotdot/.style = {
        regularline,
        dash pattern=on 3pt off 1.4pt on 0.7pt off 1.4pt on 0.7pt off 1.4pt,
        line cap=round
    },
    plotline loosely dotted/.style = {
        regularline,
        dash pattern=on 0.7pt off 2.6pt,
        line cap=round
    },
}

\pgfplotsset{
    every tick/.style={black, thinline},
    axis line style=black,
    scaled x ticks=false,
    xticklabel style={/pgf/number format/.cd, fixed, precision=0, 1000 sep={}},
    scaled y ticks=false,
    yticklabel style={/pgf/number format/.cd, fixed, precision=0, 1000 sep={}},
    subplot/.style={
        width=3.8cm,
        height=3.8cm,
        enlarge x limits=0.1,
        enlarge y limits=0.1,
        major tick length=0.1cm,
        thinline,
        plot coordinates/math parser=false,
        trim axis left,
        scale only axis,
        title style={font=\footnotesize, yshift=-1.0ex},
        label style={font=\scriptsize},
        ticklabel style={font=\scriptsize},
        legend style={
            font=\scriptsize,
            thinline, at={($(0,0)+(1cm,1cm)$)},
            anchor=center,
            align=center,
            legend cell align={left},
            /tikz/every even column/.append style={column sep=0.5cm},
        },
    },
    success_rate_subplot/.style={
        subplot,
        xmin=0,
        xmax=60,
        xtick={0,15,...,60},
        restrict x to domain=0:60,
        xlabel={Time (Minutes)},
        xlabel shift=-0.5ex,
        ymin=0.0,
        ymax=1.0,
        ytick distance=0.2,
        yticklabel={\pgfmathparse{\tick*100}\pgfmathprintnumber{\pgfmathresult}\%},
        ylabel={Percentage Solved},
        ylabel shift=-0.8ex,
    },
}

\usepackage{booktabs}
\usepackage{bbm}
\begin{document}

\maketitle

\begin{abstract}
    Column generation and branch-and-price (B\&P) are leading mathematical optimization methods for large-scale exact optimization, iterating between solving a master problem and a pricing problem.
    Due to the difficulty of discrete optimization, high-performance column generation often relies on a custom pricing algorithm built specifically to exploit the problem's structure.
    This bespoke nature of the pricing solver makes column generation a problem-specific method and hinders the use of generic implementations across a wide range of problems.
    We show that domain-independent dynamic programming (DIDP), a model-based paradigm for dynamic programming, can be used as a generic pricing solver.
    We develop new modeling features and a solving algorithm for DIDP to achieve better performance in typical pricing problems.
    We demonstrate that in four problem classes, our implementations of B\&P, with pricing by DIDP, empirically outperform an existing automated B\&P solver and B\&P with pricing by mixed-integer programming or constraint programming.
\end{abstract}

\section{Introduction}
\label{sec:intro}

Column generation is a method for solving a linear program by looping between solving a restricted master problem with a subset of decision variables and a pricing problem to generate variables that potentially lead to a better solution, rather than enumerating all variables in advance.
Branch-and-price (B\&P), the combination of branch-and-bound and column generation to solve mixed-integer programs (MIPs), has achieved state-of-the-art performance in solving large-scale combinatorial optimization problems, such as routing and scheduling \cite{Lubbecke:2005aa}.
Highly efficient problem-specific algorithms to solve the pricing problems have been developed, typically using dynamic programming (DP) \cite{Pugliese:2013aa,Lozano:2016aa,Sadykov:2021aa}.
However,  generalizing such algorithms to new problems is not necessarily easy due to their problem-specific nature and requires high implementation effort.
While generic column generation approaches have been investigated \cite{Gamrath:2010ab,Gualandi:2013aa}, they solve pricing problems using constraint-based approaches such as MIP and constraint programming (CP), rather than DP.

In this paper, we investigate the use of domain-independent dynamic programming (DIDP), a recently proposed model-based paradigm based on DP \cite{Kuroiwa:2025aa}, to solve the pricing problems in column generation.
By using DIDP, instead of designing and implementing problem-specific pricing algorithms, we formulate a pricing problem as a declarative DP model and apply general-purpose solvers.
We aim to move column generation closer towards a declarative model-based solving paradigm while enjoying the strengths of DP.

To achieve better performance in typical pricing problems, we develop three new modeling features in a software framework for DIDP: higher-order expressions such as filter and reduce operations, set resource variables, and the fractional knapsack expression.
We also develop a new solving algorithm for DIDP inspired by successful DP pricing algorithms.
As case studies, we formulate declarative column generation models and implement B\&P built on them in four routing and scheduling problems, using the new modeling features.
The experimental results show that the new features and algorithm contribute to better performance in multiple problem classes.
Our approach empirically outperforms generic B\&P baselines: an existing automated B\&P solver and B\&P approaches using MIP or CP for pricing.

\section{Related Work}
\label{sec:lit}

High-performance column generation often depends on specialized pricing solvers
that can exploit specific substructures,
limiting truly generic column generation solvers.
We briefly review the main automated approaches.
GCG is an open-source MIP solver that reformulates a given MIP model so that column generation can be applied \cite{Gamrath:2010ab}.
Unless a problem-specific pricing algorithm is implemented by a user, GCG solves both the master and pricing problems using SCIP \cite{Achterberg:2008aa}, an open-source MIP solver.
Similarly, Coluna,\footnote{\url{https://github.com/atoptima/Coluna.jl}} a decomposition framework for mathematical optimization in Julia, also provides automatic column generation and uses a MIP solver for pricing by default.
VRPSolver is a non-commercial but proprietary B\&P code for solving vehicle routing problems~\cite{Pessoa:2020aa}, with which users can model certain problem variants within the library's limitations.
In the CP community, previous work investigated using CP to solve pricing problems, achieving success in problems such as staff scheduling, packing, and sports scheduling \cite{Gualandi:2013aa}.

Many pricing problems can be reduced to a shortest path problem with resource constraints (SPPRC) \cite{Irnich:2005aa}.
Salani, Basso, and Giuffrida~\cite{Salani03032024} proposed PathWyse, a solver library for SPPRC, but its modeling interface is designed to handle only typical SPPRC constraints, such as capacities and time windows, and C++ customization is required for more complex variants.
The Boost Graph Library,\footnote{\url{https://www.boost.org/doc/libs/1_91_0/libs/graph/doc/index.html}} a C++ library for graph algorithms and data structures, also implements a function to solve SPPRC.
However, it is less declarative than model-based paradigms such as MIP, CP, and DIDP because a problem is defined by implementing a set of functions in C++.

\section{Background}
\label{sec:background}

We formally introduce column generation and domain-independent dynamic programming.

\subsection{Column Generation}

Many combinatorial optimization problems can be posed as a set partitioning or set covering formulation that selects
a set of combinatorial objects (e.g., a path, schedule, cutting pattern) from a large universe.
MIP formulations of these problems often associate each object with a variable/column.
The linear relaxation contains an exponential number of columns, so it cannot be constructed explicitly.
Instead, column generation solves a sequence of smaller linear relaxations over a subset of columns,
and dynamically generates new columns \cite{Lubbecke:2005aa}.

Let $X$ be the set of all columns, with $n = |X|$. Define the \emph{master problem (MP)} as
$\min \{ c^\top \lambda : A\lambda \geq b, \lambda \geq 0 \}$,
where $\lambda = (\lambda_1, \ldots, \lambda_n)$, $c \in \mathbb{Q}^n$, $A \in \mathbb{Q}^{m \times n}$, and $b \in \mathbb{Q}^{m}$.
Column generation maintains a subset $X' \subseteq X$ with $|X'| = n' \leq n$ and solves the corresponding
\emph{restricted master problem (RMP)}
$\min \{ c'^\top \lambda' : A'\lambda' \geq b, \lambda' \geq 0 \}$,
where $c' \in \mathbb{Q}^{n'}$ and $A' \in \mathbb{Q}^{m \times n'}$ are the corresponding submatrices.
Let $\hat{\pi} \in \mathbb{R}^m_+$ be an optimal dual solution to the RMP.
For any column $j \in X$, define its reduced cost as
$\bar c_j = c_j - A_{\cdot,j}^\top \hat{\pi}$.
Any column $j \in X \setminus X'$ with $\bar c_j < 0$
may be needed in a lower-cost solution to the RMP,
and therefore should be added to $X'$.

Since the full set of columns is too large to enumerate,
columns with negative reduced cost are found by solving a \emph{pricing problem}.
Assume that each column is represented by a vector $x \in \mathbb{Z}^k$ satisfying
internal constraints $Dx \geq e$, with cost $c(x)$ and column coefficients $a(x)$.
Then, the pricing problem is
$\min \{ c(x) - a(x)^\top \hat{\pi} : Dx \geq e, x \in \mathbb{Z}^k \}$.
If the minimum is negative, the corresponding column is added to $X'$, and the process repeats.
Otherwise, the current optimal RMP solution is also optimal for the MP.
Branch-and-price (B\&P) embeds column generation inside branch-and-bound to enforce integrality constraints
on MP.
%A full review of branching rules is beyond the scope of this work.

\subsection{Dynamic Programming}

Combinatorial problems can be defined in DP by states and transitions, with costs or profits on transitions.
As a running example, consider SPPRC instantiated as a pricing problem for the vehicle routing problem with time windows (VRPTW) \cite{Vigo:2014aa} with capacity and time-window resources and an elementary (no-revisit) constraint. Let $(\mathcal{N}, \mathcal{A})$ be a directed graph with nodes $\mathcal{N}=\{0,\ldots,n+1\}$ (source $0$, sink $n+1$) and arcs $\mathcal{A} \subseteq \{ \mathcal{N} \times \mathcal{N} : i \neq j, i < n+1, j > 0 \}$. Each customer $i$ has load $l_i > 0$, release time $a_i \ge 0$, deadline $b_i \ge 0$, and service duration $s_i \ge 0$. Each arc $(i,j)$ has distance $d_{i,j} \ge 0$ and travel cost $c_{i,j}$. The goal is an elementary path from $0$ to $n+1$ of minimum total travel cost such that cumulative load never exceeds capacity $Q$ and the visit to each node $i$ occurs within $[a_i,b_i]$.

Let $V(\mathcal{R}, i, q, t)$ be the minimum cost from node $i$ to $n+1$ when the unvisited set is $\mathcal{R}\subseteq\{1,\ldots,n\}$, current load is $q$, and time is $t$.
We may visit node $j$ only if $j \in \mathcal{R} \cup \{ n + 1 \}$, $(i, j) \in \mathcal{A}$, $q + l_j \leq Q$, and $t + s_i + d_{i,j} \leq b_j$.
Once visiting $j$, we consider a subproblem of computing the minimum cost from node $j$ to $n + 1$, where $j$ is removed from the unvisited set, the load is increased by $l_j$, and the time is updated to the later of the arrival time $t + s_i + d_{i,j}$ and the beginning of the time window $a_j$ at $j$.
As a result, we compute $V(\mathcal{R} \setminus \{ j \}, j, q + l_j, \max\{ t + s_i + d_{i,j}, a_j \})$.
Using the set of nodes that can be visited next, $\text{Next}(\mathcal{R}, i, q, t) = \{ j \in \mathcal{R} \cup \{ n + 1 \} : (i, j) \in \mathcal{A} \land q + l_j \leq Q \land t + s_i + d_{i,j} \leq b_j \}$, $V$ is defined by the following Bellman equation:
\begin{equation}
    V(\mathcal{R}, i, q, t) = \begin{cases}
            \min\limits_{j \in \text{Next}(\mathcal{R}, i, q, t) } c_{i, j} + V \left( \mathcal{R} \setminus \{ j \}, j, q + l_j, \max \left\{ t + s_i + d_{i,j}, a_j \right\} \right) & \text{if } i \leq n \\
            0 & \text{else.}
    \end{cases}
    \label{eq:example_bellman}
\end{equation}
The optimal objective value is $V(\{ 1, \cdots, n \}, 0, 0, 0)$.
In column generation, Bellman equations are typically solved using problem-specific algorithms \cite{Pugliese:2013aa,
Lozano:2016aa,Sadykov:2021aa}.

\subsection{Domain-Independent Dynamic Programming}

Domain-independent dynamic programming (DIDP) is a model-based DP paradigm for combinatorial optimization.
Previous work developed Dynamic Programming Description Language (DyPDL)~\cite{Kuroiwa:2025aa}, a declarative modeling formalism for DIDP.
Currently, DyPDL is limited to DP formulations where a solution can be represented as a sequence of decisions, and a DP model is defined by \emph{state variables}, \emph{transitions}, \emph{base cases}, and \emph{state constraints}.

A state variable has a type, either numeric, element, or set.
A \emph{numeric variable} takes a value in $\mathbb{Q}$, an \emph{element variable} in $\mathbb{Z}^+_0$, and a \emph{set variable} in $2^{\mathbb{Z}^+_0}$.
A state is represented by full value assignments $(d_1, \ldots, d_n)$ to the state variables $(x_1, \ldots, x_n)$ such that $d_i$ is in the domain of $x_i$, and we denote the value of a state variable $x_i$ in state $S$ by $S[x_i]$.
For our example in \Cref{eq:example_bellman}, $\mathcal{R}$, $i$, $q$, and $t$ are modeled as state variables, where $\mathcal{R}$ is a set variable, $i$ is an element variable, and $q$ and $t$ are numeric variables.
%Solving a DIDP model is to compute the value of a single particular state called the \emph{target state}.
%In our example, the target state is defined by $(\mathcal{R} = \{ 1, ..., n \}, i = 0, q = 0, t = 0)$.

\emph{Expressions} are functions of a state, where \emph{numeric expressions} return a numeric value, \emph{element expressions} return a nonnegative integer, \emph{set expressions} return a set of nonnegative integers, and \emph{conditions} return a Boolean value ($\bot$ or $\top$).
Given a condition $c$, we denote $S \models c$ if $c(S) = \top$, i.e., $S$ satisfies $c$, and $S \not \models c$ if $c(S) = \bot$.
In practice, expressions are built from a predefined set of operations such as arithmetic and set-theoretic operations.

A transition defines \emph{effects} on the state variables and \emph{preconditions} to apply, using expressions.
In our example, we define the first line of \Cref{eq:example_bellman} by using a transition corresponding to visiting each customer $j$, whose effect on the set variable $\mathcal{R}$ is $\mathcal{R} \setminus \{ j \}$ (a set expression), effect on the element variable $i$ is $j$ (an element expression), effect on the numeric variable $q$ is $q + l_j$ (a numeric expression), and effect on the numeric variable $t$ is $\max\{ t + s_i + d_{i,j}, a_j \}$ (a numeric expression).
By $S[\![\tau]\!]$, we denote the \emph{successor state} of state $S$, where the values of state variables are updated by the effects of transition $\tau$.
To visit customer $j$, it must satisfy $j \in \text{Next}(\mathcal{R}, i, q, t)$, i.e., conditions $j \in \mathcal{R} \cup \{ n + 1 \}$, $(i, j) \in \mathcal{A}$, $q + l_j \leq Q$, and $t + s_i + d_{i,j} \leq b_j$ are preconditions.
By $\mathcal{T}(S)$, we denote the set of \emph{applicable} transitions whose preconditions are satisfied by $S$.
%We assume the additive cost structure, where a numeric expression $w_\tau$ is associated with each transition $\tau$, and the cost becomes $w_\tau(S)$ plus the cost of the successor state.
%For our example, $w_\tau(\mathcal{R}, i, q, t) = d_{i,j}$.

A base case consists of a set of conditions when no further transition is applicable.
A state satisfying a base case is called a \emph{base state}.
In our example, we define the second line of \Cref{eq:example_bellman} using a base case, where the set of conditions is $\{ i = n + 1 \}$.
An \emph{$S$-solution} is a sequence of transitions $\sigma = \langle \sigma_1, ..., \sigma_m \rangle$ inducing a sequence of states $\langle S^0, ..., S^m \rangle$ such that $S^0 = S$,  $S^m$ is a base state, and $S^i = S^{i-1}[\![\sigma_i]\!]$ with $\sigma_i \in \mathcal{T}(S^{i-1})$ for $i = 1, \ldots, m$.
We assume the additive cost structure: the cost of the $S$-solution $\sigma$ is $\sum_{i = 1}^{m} w_{\sigma_i}(S^{i-1}) + w_\text{base}(S^m)$, where $w_\tau$ is a numeric expression associated with transition $\tau$, and $w_\text{base}$ is a numeric expression for the base cases.
In our example, $w_\tau(\mathcal{R}, i, q, t) = c_{i,j}$ for a transition $\tau$ to visit location $j$, and $w_\text{base}(\mathcal{R}, i, q, t) = 0$.
Assuming minimization, by $V(S)$, we denote the minimum $S$-solution cost.
A special state called the \emph{target state} $S^I$ is defined for a DyPDL model, and the optimal solution for a DyPDL model is the minimum cost $S^I$-solution.
In our example, $S^I = (\mathcal{R} = \mathcal{N}, i=0, q=0, t=0)$.
In addition to the above components, \emph{state constraints} are conditions that must be satisfied by all states, but they do not appear in our example.

\subsubsection{Redundant Information}

In DyPDL, redundant information implied by other parts of the model can be explicitly defined, which can potentially be useful for a solver.
DyPDL provides two specific features solely for redundant information: \emph{resource variables} and \emph{dual bound functions}.

In problem-specific DP algorithms, state dominance is sometimes exploited, where one state is known to be superior to another.
For our example in \Cref{eq:example_bellman}, a state $\left( \mathcal{R}, i, q_1, t_1 \right)$ leads to a better or equal solution than a state $\left( \mathcal{R}, i, q_2, t_2 \right)$ with the same $\mathcal{R}$ and $i$ if $q_1 \leq q_2$ and $t_1 \leq t_2$.
This state dominance can be written as the following inequality:
\begin{equation}
    V(\mathcal{R}, i, q_1, t_1) \leq V(\mathcal{R}, i, q_2, t_2) \text{ if } q_1 \leq q_2 \land t_1 \leq t_2.
\end{equation}
In DyPDL, to represent state dominance, a numeric or element variable can be declared as a resource variable with a preference for less or greater.
A state $S$ is preferred over another state $S'$ if $S[r] \leq S'[r]$ for each resource variable $r$ that prefers less, $S[r] \geq S'[r]$ for each resource variable $r$ that prefers greater, and $S[x] = S'[x]$ for each non-resource variable $x$.
In our example, $q$ and $t$ are resource variables where less is preferred.

A dual bound function $\eta$ returns a lower bound $\eta(S)$ on the minimum $S$-solution cost $V(S)$ and an upper bound for maximization, similar to a dual bound in mathematical optimization.
In DyPDL, a dual bound function is described by an expression, similar to other components.
For our example in \Cref{eq:example_bellman}, since the travel cost of an arc can be negative, we can use a dual bound function that only considers the most negative incoming arc for each node $j = 1, \ldots, n$.
We can underestimate the incoming arc cost by $c^{\text{in}}_j =  \min_{(k, j) \in \mathcal{A}} c_{k,j}$ for node $j$.
We have the following dual bound function:
\begin{equation}
    V(\mathcal{R}, i, q, t) \geq \begin{cases}
        0 & \text{ if } i = n +1 \\
        \sum_{j \in \mathcal{R} } \min\left\{ c^{\text{in}}_j, 0 \right\} + c^\text{in}_{n+1} & \text{ else.}
    \end{cases}
    \label[example]{eq:example_bellman:bound}
\end{equation}
Similarly, we can use another dual bound function using $c^\text{out}_j = \min_{(j,k) \in \mathcal{A}} c_{j,k}$ instead of $c^\text{in}_j$ for $j \in \mathcal{R}$ and $c^\text{out}_i$ instead of $c^\text{in}_{n+1}$.

\subsubsection{Solvers}

Once a DP model is formulated, it is solved by general-purpose solvers.
Previous work \cite{Kuroiwa:2025aa} developed solvers based on heuristic state-space search algorithms \cite{Russel:2020a}, such as CAASDy, which uses A* \cite{Hart:1968aa}, and complete anytime beam search (CABS) \cite{Zhang:1998aa}.
In these solvers, a DP model is solved by finding a shortest path from the target state to a base state in a state space graph, where the nodes are states, and an edge $(S, S[\![\tau]\!])$ exists if $\tau \in \mathcal{T}(S)$.
Each edge $(S, S[\![\tau]\!])$ is labeled with transition $\tau$ and having the weight $w_\tau(S)$.

\section{New Modeling Features in DIDP}
\label{sec:novelty}

We extend didp-rs,\footnote{\url{https://github.com/domain-independent-dp/didp-rs}} a software implementation of DyPDL, to efficiently model and solve problems commonly seen in pricing.
We add three new features: higher-order expressions, set resource variables, and fractional knapsack expressions.

\subsection{Higher-Order Expressions}

For the SPPRC example in \Cref{eq:example_bellman}, a state is represented by a set of unvisited nodes $\mathcal{R}$, the current node $i$, the current load $q$, and the current time $t$.
While we update $\mathcal{R}$ to $\mathcal{R} \setminus \{ j \}$ when $j$ is visited, we can also remove a node $k \in \mathcal{R}$ that can no longer be visited by its deadline $b_k$.
Let $d^*_{j,k}$ be the shortest travel time from node $j$ to node $k$, which can be precomputed.
If $\max\{ t + s_i + d_{i,j}, a_j \} + s_j + d^*_{j,k} > b_k$, then node $k$ cannot be visited after visiting $j$ from the current state.
In addition, $k$ cannot be visited after $j$ if it results in overload, i.e., $q + l_j + l_k > Q$.
Thus, the update on $\mathcal{R}$ can be represented by an expression $\mathcal{R}'(j) = \left\{ k \in \mathcal{R} \setminus \{ j \} :  \max\{ t + s_i + d_{i,j}, a_j \} + s_j + d^*_{j,k} \leq b_k \land q + l_j + l_k \leq Q \right\}$.

In didp-rs, expressions are built from a library of operations on state variables.
The solvers maintain expression tree data structures and evaluate them during solving.
For set expressions, set operations such as union, intersection, and difference are implemented, but $\mathcal{R}'(j)$ cannot be directly represented with them.
%In addition, an `if-then-else' operation is available, which evaluates to one of two expressions depending on the evaluation result of a condition.
%Using these operations, $\mathcal{R}'(j)$ can be implemented by repeatedly removing a singleton or empty set defined by a set expression `if $k \in \mathcal{R} \land (t'(j) + s_j + d^*_{j,k} > b_k \lor q + l_j + l_k > Q)$ then $\{ k \}$ else $\emptyset$' for each $k = 1, \dots, n$.
%However, such an implementation complicates the code and results in an expression tree whose depth is proportional to $n$, which is slow to evaluate in practice.
%We show an example of such implementation using DIDPPy, the Python interface of didp-rs, in Listing~\ref{list:manual-filter}, where $\{ 0 \}$ is used instead of $\emptyset$ as the source node $0$ is not included in $\mathcal{R}$.
%
%\begin{lstlisting}[caption={DIDPPy example of the filter operation with existing features.},label=list:manual-filter,float=ht]
%q_new = q + l[j]
%t_new = dp.max(t + s[i] + d[i,j], a[j])
%r_new = r.remove(j)
%
%for k in range(1, n+1):
%    c = (t_new + s[j] + d_star[j,k] <= b[k]) & (q_new + l[k] <= Q)
%    r_new = r_new.remove(c.if_then_else(k, 0))
%\end{lstlisting}

We introduce higher-order expressions, which apply an expression to each element of a given set.
In particular, a \emph{filter operation} is a set expression that returns a subset of a given set whose elements satisfy a given condition.
With our interface, a user specifies a filter operation by two components: a set expression $\mathcal{X}$ and a parameterized condition $c$, which is a function that returns a condition $c(x)$ given a parameter $x$.
The parameter $x$ is a placeholder and is replaced with each element of a set $\mathcal{X}(S)$ when evaluated, given a state $S$, and an element $i \in \mathcal{X}(S)$ is removed if $S \not \models c(i)$.
In other words, the filter operation represents an expression that returns $\{ x \in \mathcal{X}(S) : S \models c(x) \}$ given a state $S$.
For our example, $\mathcal{R}'(j)$ can be represented by a filter operation defined by a set expression $\mathcal{R} \setminus \{ j \}$ and a parameterized condition $\max\{ t + s_i + d_{i,j}, a_j \} + s_j + d^*_{j,k} \leq b_k \land q + l_j + l_k \leq Q$, where $k$ is the parameter.
We present example code using the filter operation in DIDPPy, the Python interface of didp-rs, in Listing~\ref{list:filter}.
Here, we define the placeholder using \textsf{model.add\_local\_var} and give it to a filter expression (\textsf{.filter}) as the first argument, where \textsf{r} is a set variable.

\begin{lstlisting}[caption={DIDPPy example of the filter operation.},label=list:filter,float=ht]
q_next = q + l[j]
t_next = dp.max(t + s[i] + d[i,j], a[j])
k = model.add_local_var()
r_next = r.remove(j).filter(
    k, (t_next + s[j] + d_star[j,k] <= b[k]) & (q_next + l[k] <= Q)
)
\end{lstlisting}

We also introduce the \emph{reduce sum operation} returning $\sum_{x \in \mathcal{X}(S)} e(x)(S)$ given a state $S$, where $\mathcal{X}$ is a set expression, $x$ is a placeholder, and $e$ is a parameterized numeric expression.
We show its use case later in Section~\ref{sec:mrasp}.

\subsection{Set Resource Variables}

In the current DyPDL, only numeric and element variables can be resource variables to define state dominance.
However, in pricing problems, state dominance is sometimes defined by a set variable.
For the SPPRC example in \Cref{eq:example_bellman}, we can define state dominance where state $(\mathcal{R}_1, i, q_1, t_1)$ is better than or as good as $(\mathcal{R}_2, i, q_2, t_2)$ with the same $i$ if $\mathcal{R}_2 \subseteq \mathcal{R}_1$, $q_1 \leq q_2$, and $t_1 \leq t_2$ since having more candidates to visit potentially leads to a shorter path.

We introduce \emph{set resource variables}: a state $S$ is preferred to another state $S'$ only if the value of a set resource variable in $S$ is a subset or superset of that in $S'$.
Similar to numeric and element resource variables, the preference (less or greater) specifies whether a subset or superset is better.
When less/greater is specified for a set resource variable $\mathcal{X}$, $S$ is preferred to $S'$ only if $S[\mathcal{X}] \subseteq S'[\mathcal{X}]$/$S'[\mathcal{X}] \subseteq S[\mathcal{X}]$.
%A set resource variable can be mimicked by defining a set of numeric or element resource variables, whose values take either $0$ or $1$.
%However, our set resource variable implementation uses a bitset to represent a set, which is computationally more efficient.
We present example code in Listing~\ref{list:set-resource}, where \textsf{node} is an object type to define the maximum cardinality of the set variable, \textsf{target} specifies the value in the target state, and \textsf{less\_is\_better=False} specifies that greater is preferred.

\begin{lstlisting}[caption={DIDPPy example of a set resource variable},label=list:set-resource,float=ht]
r = model.add_set_resource_var(
    node, target=list(range(1, n + 1)), less_is_better=False
)
\end{lstlisting}

\subsection{Fractional Knapsack Expression}

For our example in \Cref{eq:example_bellman}, we presented  in \Cref{eq:example_bellman:bound} a dual bound function based on the minimum incoming arc cost $c^{\text{in}}_j = \min_{(k,j) \in \mathcal{A}} c_{k,j}$ for each node $j$.
We can also define a dual bound based on the current load $q$ and the capacity $Q$.
By visiting node $j$, we increase the load by $l_j$ and the cost by at least $c^{\text{in}}_j$.
Given a state $(\mathcal{R}, i, q, t)$,
\begin{align}
    V(\mathcal{R}, i, q, t) \geq  \min_{\mathcal{J} \subseteq \mathcal{R} : q + \sum_{j \in \mathcal{J}} l_j \leq Q} \sum_{j \in \mathcal{J}} \min\left\{ c^{\text{in}}_j, 0 \right\} + c^\text{in}_{n+1} & \text{ if } i \neq n +1.
    \label[ineq]{eq:example_bellman:knapsack_bound}
\end{align}
The first term can be viewed as the negation of the optimal cost of the 0-1 knapsack problem, which is to maximize the total profit of items packed into a knapsack with a fixed capacity.
In particular,  the knapsack has capacity $Q - q$, and each node $j \in \mathcal{R}$ with $c^{\text{in}}_j < 0$ corresponds to an item with the profit $-c^{\text{in}}_j$ and weight $l_j$.
We argue that a similar substructure is common in pricing problems when a subset of elements with the negative reduced costs needs to be selected under a resource constraint.

Since the 0-1 knapsack problem is NP-hard~\cite{Karp:1972aa}, computing the right-hand side of \Cref{eq:example_bellman:knapsack_bound} is also NP-hard.
We use the Dantzig bound \cite{Dantzig:1957aa}, a polynomial-time upper bound on the optimal objective value for the 0-1 knapsack problem.
%Recent work has reported that the Dantzig bound~\cite{Dantzig:1957aa}, a polynomial-time upper bound on the optimal objective value for the 0-1 knapsack problem, is useful as the dual bound function for DIDP using RPID,  \cite{Kuroiwa:2025ab}.
Given the capacity $C$ and a set of items $\mathcal{N}$ with weight $w_j > 0$ and the profit $p_j > 0$ for each $j \in \mathcal{N}$, the Dantzig bound can be computed as follows.
First, the items are sorted in a descending order of $\frac{p_j}{w_j}$.
Second, the items are included in the knapsack in sorted order as long as the total weight does not exceed the capacity $C$, and let $\mathcal{I}$ be the set of such items.
When the current item $j$ has the weight $w_j$ larger than $C - \sum_{i \in \mathcal{I}} w_i$, it is fractionally included with the profit $\frac{p_j}{w_j} \left( C - \sum_{i \in \mathcal{I}} w_i \right)$, i.e., the objective value is upper bounded by $\frac{p_j}{w_j} \left( C - \sum_{i \in \mathcal{I}} w_i \right) + \sum_{i \in \mathcal{I}} p_i$.

The Dantzig bound was used with RPID, another DIDP software where a model is defined by Rust functions \cite{Kuroiwa:2025ab}.
However, with expressions in didp-rs, efficiently modeling the Dantzig bound is difficult due to its algorithmic nature.
We introduce a new expression, called the fractional knapsack expression, denoted by $\textsf{fractional\_knapsack} \left( \mathcal{X}, C, (p_j)_{j = 1, \dots, n}, (w_j)_{j = 1, \dots, n} \right)$, where $\mathcal{X}$ is a set expression, $C$ is a numeric expression, and $(p_j)_{j = 1, \dots, n}$ and $(w_j)_{j=1, \dots, n}$ are lists of $n$ numeric expressions.
Then, the expression represents the Dantzig bound for the 0-1 knapsack problem, where given a state $S$, the set of items is $\mathcal{X}(S)$, the capacity of the knapsack is $C(S)$, and each item $x \in \mathcal{X}(S)$ has the profit $p_x$ and the weight $w_x$.
For our example, when $i \neq n + 1$, we represent the dual bound function as follows:
\begin{equation}
    V(\mathcal{R}, i, q, t) \geq -\textsf{fractional\_knapsack} \left( \mathcal{R}, Q - q, \left( \max \left\{ -c^{\text{in}}_j, 0 \right\} \right)_{j = 1, \dots, n}, \left( l_j \right)_{j = 1, \dots, n} \right) + c^\text{in}_{n+1}
    \label[example]{eq:example_bellman:fractional_knapsack_bound}
\end{equation}
where $c^\text{in}_j = \min_{(k,j) \in \mathcal{A}} c_{k,j}$ is the minimum incoming travel cost to node $j$.
Zero-weight items are ignored by the expression.
In our example code (Listing~\ref{list:fractional-knapsack}), \textsf{(i == n+1).if\_then\_else(..., ...)} is an expression that evaluates to the first argument (0) if the current location \textsf{i} equals \textsf{n + 1} (the sink node) and to the second argument (the fractional knapsack bound) if not.
\begin{lstlisting}[caption={DIDPPy example of a fractional knapsack expression},label=list:fractional-knapsack,float=ht]
model.add_dual_bound(
    (i == n + 1).if_then_else(
        0,
        -dp.fractional_knapsack(
            r, Q - q, [max(-cin[j], 0) for j in range(1, n + 1)], l
        )
    )
)
\end{lstlisting}

We can use three additional dual bound functions similar to Inequality~\eqref{eq:example_bellman:fractional_knapsack_bound}:
the first one uses the capacity $b_{n+1} - t -s_i$ and the weight $\min_{(k,j) \in \mathcal{A}} d_{k,j} + s_j$;
the second one uses the profit $\max\left\{ -c^\text{out}_j, 0 \right\}$ with $c^\text{out}_j = \min_{(j,k) \in \mathcal{A}} c_{j,k}$ instead of $c^\text{in}_j$;
the third one uses the profit $\max\left\{-c^\text{out}_j, 0 \right\}$, the capacity $b_{n+1} - t$, and the weight $\min_{(j,k) \in \mathcal{A}} s_j + d_{j,k}$.

\section{Generic Labeling Solver for DIDP}

Our solver is built on top of the anytime heuristic search framework of DIDP proposed in previous work~\cite{Kuroiwa:2025aa}, which solves a DP model by finding a shortest path from the target state to a base state in a state space graph.
In that framework, states to be searched are maintained in a priority queue called an \emph{open list}.
In each iteration, one state is selected and removed from the open list.
Then, it is \emph{expanded}, i.e., its successor states are generated by applying transitions and then inserted into the open list if they are not dominated by existing states.
Each concrete algorithm differs in how it selects states to expand from the open list.
In existing DIDP solvers, including CAASDy and CABS, the state is selected based on its $f$-value, which is the sum of its $g$-value, the path cost to reach the state from the target state, and its $h$-value, the estimated path cost from the state to a base state by a heuristic function $h$.
The existing solvers use the dual bound function $\eta$ defined in the model as the heuristic function.
When multiple dual bound functions are defined, the maximum value is used for each state.

Unlike the existing solvers in DIDP, our new solver is inspired by problem-specific labeling algorithms for solving the SPPRC~\cite{Irnich:2005aa,Pugliese:2013aa} and selects states based on resource variables.
We explain the motivation of our labeling solver using the example SPPRC for VRPTW.
Suppose two states $S_1 = (\mathcal{R}_1, i, q_1, t_1)$ with $g$-value $g(S_1)$ and $S_2 = (\mathcal{R}_2, i, q_2, t_2)$ with $g$-value $g(S_2)$ such that $\mathcal{R}_2 \subseteq \mathcal{R}_1$, $q_2 \geq q_1$, $t_2 \geq t_1$, and $g(S_2) \geq g(S_1)$, i.e., $S_1$ dominates $S_2$ with a better or equal $g$-value.
Since extending the path from the target state to $S_1$ results in a better or equal solution than that of $S_2$, we do not need to search $S_2$ once expanding $S_1$.
In other words, expanding $S_2$ before $S_1$ should be avoided, and expanding the state with more preferred resource variable values potentially reduces search effort.
With this motivation, our algorithm expands a state based on a lexicographic order of resource variables.

\begin{algorithm}[t]
    \caption{
        Generic labeling solver for a DyPDL model.
        The target state is denoted by $S^I$ and the dual bound function by $\eta$.
    }
    \begin{algorithmic}[1]
        \If{$S^I$ does not satisfy the state constraints}
        \Return $\emptyset$ \label{alg:search:return-empty} %\Comment{Check the state constraints.}
        \EndIf
        \State $\Sigma \leftarrow \emptyset$, $\overline{\gamma} \leftarrow \infty$, $\sigma(S^I) \leftarrow \langle \rangle$, $g(S^I) \leftarrow 0$, $O \leftarrow \{ S^I \}, G \leftarrow \{ S^I \}$ \label{alg:search:init} \Comment{Initialization.}
        \While{$O \neq \emptyset$} \label{alg:search:while}
        \State Let $S \in O$ be the lexicographically minimum state \label{alg:search:select}
        \State $O \leftarrow O \setminus \{ S \}$ \label{alg:search:remove} \Comment{Pop the chosen state.}
        \If{$S$ is a base state and $g(S) + w_\text{base}(S) < \overline{\gamma}$} \label{alg:search:check-base}
        \State $\Sigma \leftarrow \Sigma \cup \{ \sigma(S) \}$, $\overline{\gamma} \leftarrow g(S) + w_\text{base}(S)$ \label{alg:search:update-solution} \Comment{Update the solutions.}
        \State $O \leftarrow \{ S' \in O : g(S') + \eta(S') < \overline{\gamma} \}$ \label{alg:search:prune-open} \Comment{Prune states in the open list.}
        \Else
        \ForAll{$\tau \in \mathcal{T}(S) : S[\![\tau]\!]$ satisfies all state constraints} \label{alg:search:gen}
        \State $g_{\text{current}} \leftarrow g(S) + w_\tau(S)$ \Comment{Compute the $g$-value.}
        \If{no state $S' \in G$ is preferred to $S[\![\tau]\!]$ with $g(S') \leq g_\text{current}$} \label{alg:search:dominating}
        \State $G \leftarrow \{ S' \in G : S[\![\tau]\!] \text{ is not preferred to } S' \lor g(S') < g_\text{current}  \}$ \label{alg:search:dominated}
        \State $O \leftarrow \{ S' \in O : S[\![\tau]\!] \text{ is not preferred to } S' \lor g(S') < g_\text{current}  \}$ \label{alg:search:dominated-open}
        \If{$g_{\text{current}} + \eta(S[\![\tau]\!]) < \overline{\gamma}$} \label{alg:search:check-bound}
        \State $\sigma(S[\![\tau]\!]) \leftarrow \langle \sigma(S); \tau \rangle$, $g(S[\![\tau]\!]) \leftarrow g_{\text{current}}$ \label{alg:search:update-g}
        \State $G \leftarrow G \cup \{ S[\![\tau]\!] \}$, $O \leftarrow O \cup \{ S[\![\tau]\!] \}$ \label{alg:search:insert} \Comment{Insert the successor state.}
        \EndIf
        \EndIf
        \EndFor
        \EndIf
        \EndWhile
        \State \Return $\Sigma$ \label{alg:search:return} \Comment{Return solutions.}
    \end{algorithmic}
    \label{alg:labeling}
\end{algorithm}

We present pseudocode in Algorithm~\ref{alg:labeling}.
Except for line~\ref{alg:search:select}, the algorithm and implementation details mirror the existing solvers in DIDP.
When the target state violates state constraints, we immediately return an empty set and terminate (line~\ref{alg:search:return-empty}).
%In column generation, generating multiple columns in one iteration of pricing is typically beneficial.
%To emphasize that the algorithm returns multiple solutions,
We initialize the set of solutions found $\Sigma$ with an empty set (line~\ref{alg:search:init}).
We also maintain the current best solution cost $\overline{\gamma}$, initialized with $\infty$ (line~\ref{alg:search:init}).
For each state $S$, we record the best sequence of transitions to reach it, $\sigma(S)$, and the $g$-value $g(S)$, corresponding to the accumulated path cost.
Given $\sigma(S) = \langle \sigma_1, \ldots, \sigma_m \rangle$, we have $g(S) = \sum_{i=1}^m w_{\sigma_i}(S^{i-1})$ where $S^i = S^{i-1}[\![\sigma_i]\!]$ for $i = 1, \ldots, m$ and $S^0 = S^I$.
For the target state $S^I$, the path is empty with $g(S^I) = 0$.
The set $G$ stores all generated states to check state dominance with a newly generated state, and the open list $O$ stores states to be searched, both of which initially contain only the target state.
The algorithm proves optimality (or infeasibility) when the open list becomes empty (line~\ref{alg:search:while}) and returns the set of solutions found.

In each step, the lexicographical-minimum state $S$ is removed from the open list (lines~\ref{alg:search:select} and \ref{alg:search:remove}).
States are lexicographically ordered based on the values of resource variables.
Given resource variables $r_1, \ldots, r_{n'}$, a state $S$ is lexicographically smaller than $S'$ if there exists $1 \leq i \leq n'$ such that $S[r_j] = S'[r_j]$ for $1 \leq j < i$, $S[r_i] \neq S'[r_i]$, and $S[r_i]$ is preferred to $S'[r_i]$.
In our implementation, we compare element resource variables, numeric resource variables, and set resource variables in order.
Resource variables of the same type are compared in order of definition.
When all resource variables have the same values, we break ties by the $f$-value ($f(S) = g(S) + \eta(S)$), and then the $\eta$-value, where smaller is preferred.
When no resource variables are used, our algorithm is the same as CAASDy, which uses A*.

If $S$ is a base state, then $\sigma(S)$ is a solution, and the best solution cost $\overline{\gamma}$ is updated if $\sigma(S)$ is better (lines~\ref{alg:search:check-base}--\ref{alg:search:update-solution}).
In addition, all states $S' \in O$ with $g(S') + \eta(S') \geq \overline{\gamma}$ are removed from the open list since they cannot lead to a better solution (line~\ref{alg:search:prune-open}).
%Here, $g(S') + \eta(S')$ can be viewed as the $f$-value $g(S') + h(S')$, but we use it only for pruning states.
%Since $\eta(S)$ is a lower bound on the solution cost starting from $S$, $g(S) + \eta(S)$ is a lower bound on the solution cost extending the sequence of transitions $\sigma(S)$.
%If this value is equal to or worse than the current solution cost, the current sequence does not lead to a better solution, so we ignore it.

If $S$ is not a base state, its successor state $S[\![\tau]\!]$ is generated for every applicable transition in $\tau \in \mathcal{T}(S)$ if it satisfies the state constraints (line~\ref{alg:search:gen}).
If $S[\![\tau]\!]$ is dominated by another state $S'$ in $G$ with a better or equal $g$-value, it cannot lead to a solution better than $S'$, so $S[\![\tau]\!]$ is ignored (line~\ref{alg:search:dominating}).
Otherwise, states dominated by $S[\![\tau]\!]$ with a better or equal $g$-value are removed from $G$ (line~\ref{alg:search:dominated}).
For this procedure, $G$ is implemented as a hash table, where keys are the values of the non-resource variables, and entries are arrays of pointers to states.
When a successor state is generated, an array of states with the same non-resource variable values is retrieved from the hash table.
The successor state is compared against each state in the array to detect dominance and appended to the array if not dominated.

After dominance detection, the dual bound value $\eta(S[\![\tau]\!])$ is computed.
If $g(S) + w_\tau(S) + \eta(S[\![\tau]\!])$ is worse than the best solution cost, the successor state $S[\![\tau]\!]$ is ignored (line~\ref{alg:search:check-bound}).
Otherwise, $\sigma(S[\![\tau]\!])$ and $g(S[\![\tau]\!])$ are initialized or updated, and $S[\![\tau]\!]$ is inserted into the open list and $G$ (line~\ref{alg:search:insert}).
Here, $\langle \sigma(S); \tau \rangle$ is a sequence of transitions extending $\sigma(S)$ with $\tau$.

%\subsection{Summary of the New Features}
%
%In summary, we add the following new features:
%\begin{itemize}
%    \item The filtering operation to efficiently construct a subset of elements satisfying a given condition.
%    \item Set resource variables for dominance pruning.
%    \item The fractional knapsack expression to efficiently compute an informative dual bound.
%    \item A generic labeling solver considering resource variables in search order.
%\end{itemize}
%The filtering operation potentially reduces the size of the state space by removing unnecessary elements from set state variables.
%Together with the filtering operation, a set resource variable enables the solving algorithm to detect more state dominance.
%The fractional knapsack expression can provide an informative dual bound.
%These two features are useful for the generic labeling solver (and other solvers) to prune unnecessary states, as shown in Algorithm~\ref{alg:labeling}.
%Furthermore, the generic labeling solver tries to avoid expanding states dominated by other states generated later.
%By combining the new modeling features and the new solver, we generalize labeling algorithms used in problem-specific settings to DIDP.

\section{Case Studies} \label{sec:case-studies}

We present DIDP models for pricing across two scheduling problems.

\subsection{Parallel Machine Scheduling} \label{sec:pmwc}

We consider the parallel machine scheduling problem commonly denoted as $P||\sum w_j C_j$~\cite{Eastman:1964aa}.
In this problem, a set of $n$ non-preemptive jobs $\mathcal{J}$ are scheduled on $m$ identical machines, where each job $j \in \mathcal{J}$ has processing time $p_j$ and weight $w_j$.
Given a schedule, let $C_j$ be the completion time of job $j$.
The objective is to minimize the total weighted completion time $\sum_{j \in \mathcal{J}} w_j C_j$.
Previous work studied B\&P for this problem \cite{VanDenAkker:1999a}, which used DP to solve the pricing problem.
However, unlike the original DP algorithm that does not use the dual bound function, our DP model uses a fractional knapsack expression as the dual bound function.

\Cref{eq:cg_pwc} shows the MP.
Let $\mathcal{S}$ be the set of schedules for a single machine, $c_s$ be the cost of the schedule, $a_{s,j}$ equal to 1 if schedule $s$ contains job $j$ and 0 otherwise, and $\lambda_s$ be a binary variable representing whether schedule $s \in \mathcal{S}$ is used.
\Cref{eq:cg_pwc:objective} minimizes the total cost, \Cref{eq:cg_pwc:machine-constraint} ensures that at most $m$ schedules are used, and \Cref{eq:cg_pwc:job-constraint} ensures that each job $j \in \mathcal{J}$ is scheduled exactly once.
\begin{subequations}
    \label[prob]{eq:cg_pwc}
    \begin{flalign}
        & \min \sum_{s \in \mathcal{S}} c_s \lambda_s
        \label[target]{eq:cg_pwc:objective} \\
        & \sum_{s \in \mathcal{S}} \lambda_s \leq m
        \label[constr]{eq:cg_pwc:machine-constraint} \\
        & \sum_{s \in \mathcal{S}} a_{s,j} \lambda_s = 1 & \forall j \in \mathcal{J}
        \label[constr]{eq:cg_pwc:job-constraint} \\
        & \lambda_s \in \mathbb{Z}_+ & \forall s \in \mathcal{S}.
    \end{flalign}
\end{subequations}

The pricing problem finds a schedule that minimizes the reduced cost, computed from the dual value $\pi_j$ for each job $j$ defined by \Cref{eq:cg_pwc:job-constraint} and the dual value defined by \Cref{eq:cg_pwc:machine-constraint}.
Without loss of generality, we assume that jobs are ordered so that $w_j / p_j \geq w_{j+1} / p_{j+1}$ for $j = 1, ..., n-1$, and job $j$ is scheduled earlier than job $k > j$ on the same machine \cite{Elmaghraby:1974aa}.
In addition, we can derive the release date (the minimum start time) $a_j$ and the deadline (the maximum completion time) $b_j$ of job $j$ based on a theoretical analysis \cite{VanDenAkker:1999a}.
Therefore, we decide whether to include job $j$ in the schedule, starting from $j = 0$ and ensuring that the starting time is in time window $[a_j, b_j]$.

In our DP formulation shown in \Cref{eq:dp_pwc}, an element variable $j$ represents the job currently considered, and a numeric variable $t$ represents the current time.
The first line of \Cref{eq:dp_pwc:transition} is a base case, where all jobs have been considered.
The second line considers a situation where $j$ can be scheduled satisfying the time window constraint.
It takes the minimum of the two cases: one where $j$ is scheduled and the other where $j$ is not scheduled.
While $j$ is updated to $j+1$ in any case, the cost is increased by $-\pi_j + w_j (t + p_j)$ and $t$ is increased by $p_j$ if $j$ is scheduled.
The third line ignores job $j$ as it cannot be scheduled within the time window.
The reduced cost is computed by subtracting the dual value for \Cref{eq:cg_pwc:machine-constraint} from $V(j = 1, t = 0)$.

An upper bound on the completion time of a schedule is given by $H = \sum_{k \in \mathcal{J}} p_k / m + (m - 1) \max_{k \in \mathcal{J}} p_k /m$.
Therefore, the sum of processing time scheduled from state $(j, t)$ must be less than or equal to $H - t$.
The completion time of each job is at least $a_k + p_k$.
Given these values, in \Cref{eq:dp_pwc:bound}, we consider a dual bound function based on the 0-1 knapsack problem with the set of items $\{j, ..., n\}$ and the capacity $H - t$, where each job $k$ has the profit $v_k = \max\{ \pi_k - w_k(a_k+p_k), 0 \}$ and the weight $p_k$.
\begin{subequations}
    \label[prob]{eq:dp_pwc}
    \begin{align}
        & V(j, t) = \begin{cases}
            0 & \text{if } j = n + 1 \\
            \min \left\{- \pi_j + w_j (t + p_j) +  V(j + 1, t + p_j) , V(j + 1, t) \right\} & \text{if } a_j \leq t \leq b_j - p_j \\
            V(j + 1, t) & \text{else}
        \end{cases}
        \label{eq:dp_pwc:transition} \\
        &  V(j, t)  \geq  -\textsf{fractional\_knapsack} \left( \{ k \in \mathcal{J} : k \geq j \}, H -t,  (v_k)_{k \geq j}, (p_k)_{k \geq j} \right).
        \label[ineq]{eq:dp_pwc:bound}
    \end{align}
\end{subequations}

Following the previous work \cite{VanDenAkker:1999a}, we change the time window of one job at each branching.

\subsection{Multi-runway Aircraft Scheduling} \label{sec:mrasp}

The multi-runway aircraft scheduling problem (MRASP), described by~\cite{Ghoniem:2015aa}, schedules a set of heterogeneous aircraft on a set of identical runways while respecting minimum separation times between aircraft and minimizing a weighted sum of scheduled times.
Let $\mathcal{M} = \{ 1,\ldots,m \}$ be the set of $m$ identical runways and let $\mathcal{O} = \{ \text{takeoff}, \text{landing} \}$ be the set of operations.
Let $\mathcal{G}$ be the set of aircraft classes and $\mathcal{N} = \{ 1,\ldots,n \}$ be the set of aircraft.
Every aircraft $i \in \mathcal{N}$ is associated with a class $g_i \in \mathcal{G}$, an operation $o_i \in \mathcal{O}$, a release time $a_i \geq 0$, a deadline $b_i \geq a_i$ and a cost weight $w_i > 0$.
Every tuple $(g_i,o_i,g_j,o_j) \in \mathcal{G} \times \mathcal{O} \times \mathcal{G} \times \mathcal{O}$
is associated with a minimum separation time $d_{g_i,o_i,g_j,o_j} \geq 0$.
An aircraft $j \in \mathcal{N}$ scheduled sometime after $i \in \mathcal{N}$, $i \neq j$, on the same runway must occur at least $d_{g_i, o_i, g_j, o_j}$ later.
In general, the triangle inequality does not hold in $d$.
Therefore, it is insufficient to check the minimum separation time of every aircraft and its immediate successor.
Instead, the minimum separation time of every later aircraft must be checked for every aircraft.

Previous work proposed B\&P for MRASP \cite{Ghoniem:2015aa}.
As MRASP can be viewed as a variant of parallel machine scheduling, where runways are machines and aircraft operations are jobs, the MP is the same as \Cref{eq:cg_pwc}.
The pricing problem is to find a schedule for a single runway satisfying the separation time constraints, for which the previous work used DP as a specialized algorithm but did not provide a declarative DP model.

In our DP model in \cref{eq:dp_mrasp}, we use a set variable $\mathcal{R}$ representing the set of aircraft whose operations can be conducted, an element variable $i$ representing the last aircraft landed, and a numeric variable $t$ representing the time to start the operation for $i$.
Each transition schedules an operation for aircraft $j \in \mathcal{R}$, removing $j$ from $\mathcal{R}$ and updating $i$ to $j$.
If the triangle inequality holds for the separation time, we can start the operation for $j$ at time $\max\{ t + d_{g_i,o_i,g_j,o_j}, a_j \}$, considering only the separation time between $i$ and $j$.
However, since the triangle inequality does not hold, depending on the aircraft operations scheduled before $i$, we may need to start the operation for $j$ later.
We use additional numeric variables $(e_{g,o})_{(g,o) \in \mathcal{Q}}$ representing the earliest time that class $g$ can perform operation $o$, where
\begin{equation*}
\mathcal{Q} = \left\{ (g,o) \in \mathcal{G} \times \mathcal{O} :  
\exists g_1 \in \mathcal{G}, o_1 \in \mathcal{O}, g_2 \in \mathcal{G}, o_2 \in \mathcal{O},
d_{g_1,o_1,g_2,o_2} + d_{g_2,o_2,g,o} < d_{g_1,o_1,g,o} \right\}.
\end{equation*}
In other words, we maintain $e_{g,o}$ for $(g, o) \in \mathcal{Q}$ such that there exist pairs of aircraft classes and operations violating the triangle inequality.
Given these variables, $t$ is updated to $t'(j) = \max\{ t + d_{g_i,o_i,g_j,o_j}, a_j, e_{g_j,o_j} \}$ if $(g_j,o_j) \in \mathcal{Q}$ and $t'(j) = \max\{ t + d_{g_i,o_i,g_j,o_j}, a_j \}$ otherwise.
Each $e_{g,o}$ is updated to $e'_{g,o}(j) = \max\{ e_{g,o}, t'(j) + d_{g_j,o_j,g,o} \}$.
For $\mathcal{R}$, we exclude aircraft $k$ whose operation cannot be executed by the deadline, i.e., $t'(j, k) > b_k$, where $t'(j, k) = \max\{ t'(j) + d_{g_j,o_j,g_k,o_k}, e'_{g_k,o_k}(j) \}$ if $(g_k,o_k) \in \mathcal{Q}$ and $t'(j,k) = t'(j) + d_{g_j,o_j,g_k,o_k}$ otherwise.
In addition, if the cost of scheduling aircraft $k$ cannot be negative, i.e., $-\pi_k + w_k t'(j, k) > 0$, there is no benefit to consider $k$.
Therefore, $\mathcal{R}$ is updated to $\mathcal{R}'(j) = \{ k \in \mathcal{R} \setminus \{ j \} : t'(j, k) \leq b'_k \}$, where $b'_k = \min\{ b_k, -\pi_k / w_k \}$, implemented by the filter operation.

\begin{subequations}
    \label[prob]{eq:dp_mrasp}
    \begin{align}
        & \begin{array}{l}
        V \left( \mathcal{R}, i, t, \left( e_{g,o} \right)_{(g,o) \in \mathcal{Q}} \right) = \\
        \begin{cases}
            0 & \text{if } i = n + 1 \\
            \min\limits_{j \in \mathcal{R} \cup \{ n + 1 \} : t'(j) \leq b'_j}  - \pi_j + w_j t'(j) + V \left( \mathcal{R}'(j), j, t'(j), \left( e'_{g,o}(j) \right)_{(g,o) \in \mathcal{Q}} \right) & \text{if } i \leq n
        \end{cases}
        \end{array} \label{eq:dp_mrasp:recursive} \\
        & \begin{array}{l}
        V \left( \mathcal{R}_1, i, t_1, \left( e^1_{g,o} \right)_{(g,o) \in \mathcal{Q}} \right)
        \leq
        V \left( \mathcal{R}_2, i, t_2, \left( e^2_{g,o} \right)_{(g,o) \in \mathcal{Q}} \right) \\
        \text{if } \mathcal{R}_2 \subseteq \mathcal{R}_1 \land t_1 \leq t_2 \land e^1_{g,o} \leq e^2_{g,o}, \forall (g,o) \in \mathcal{Q}
        \end{array}
        \label[ineq]{eq:dp_mrasp:dominance} \\
        & V \left( \mathcal{R}, i, t, \left( e_{g,o} \right)_{(g,o) \in \mathcal{Q}} \right) \geq \begin{cases}
            0 & \text{ if } i = n + 1 \\
            \sum_{j \in \mathcal{R}} \min\left\{ 0, -\pi_j +  w_j t'(j) \right\} & \text{ else.}
        \end{cases}  \label[ineq]{eq:dp_mrasp:bound}
    \end{align}
\end{subequations}

We represent a transition to finish scheduling by introducing a dummy aircraft $n+1$, with $a_j = w_j = \pi_j = 0$, $b_j = \infty$, and $d_{g_i,o_i,g_{n+1},o_{n+1}} = 0$ for any $i \in \mathcal{N}$ and define a base case in the first line of \Cref{eq:dp_mrasp:recursive}.
The minimum reduced cost is computed from $V(\mathcal{R} = \{ j \in \mathcal{N} : -\pi_j + w_j a_j < 0\}, i=0, t=0, (e_{g,o}=0)_{(g,o) \in \mathcal{Q}})$, where $i=0$ is a dummy aircraft with $d_{g_0,o_0,g_j,o_j} = 0$ for any $j \in \mathcal{N}$.
We define $\mathcal{R}$ as a set resource variable where greater is preferred (\Cref{eq:dp_mrasp:dominance}).
All numeric variables are also resource variables, where less is preferred.
The dual bound function in \Cref{eq:dp_mrasp:bound} is based on the fact that the schedule time of aircraft operation $j$ is at earliest $t'(j)$.
The dual bound function takes the sum of a parameterized expression $\min\{0, -\pi_j + w_j t'(j) \}$ with a set expression $\mathcal{R}$.

Following the previous work \cite{Ghoniem:2015aa}, branching introduces constraints specifying aircraft $k$ must or must not be scheduled immediately after $j$.
With this constraint, scheduling aircraft with a positive cost can be beneficial; if aircraft $j$ with a large negative dual value cannot be scheduled after the current aircraft $i$ but can be after another aircraft $k$, we may want to schedule $k$ even if its cost is positive.
Thus, except for the root node of B\&P, we use $b'_k = b_k$ in $\mathcal{R}'(j)$ and $(\mathcal{R} = \mathcal{N}, i=0, t=0, (e_{g,o}=0)_{(g,o) \in \mathcal{Q}})$ as the target state.

\section{Empirical Evaluation}
\label{sec:results}
%We address the following two questions through an empirical evaluation:
%\begin{enumerate}
%    \item Is our approach better than other model-based B\&P approaches?
%    \item Do the new features added to DIDP improve the performance?
%\end{enumerate}
%To answer the first question, we compare B\&P with DIDP for pricing against generic B\&P baselines on four problem classes.
%For the second question, we conduct an ablation study, evaluating variants disabling one of the features.
%
%\subsection{Experimental Settings}

To evaluate the overall performance of our framework as an exact method, we apply DIDP as a pricing solver in B\&P.
We compare it with generic B\&P baselines and conduct an ablation study to evaluate the importance of the new DIDP features.
We use the following four problem classes:
\begin{itemize}
    \item \textbf{VRPTW:}
        The pricing problem is the SPPRC used in our running example.
        We present the MP in Appendix~\ref{sec:additional-models}.
        We use the Solomon~\cite{Solomon:1987aa} instances with 50 and 100 customers. There are 112 instances in total.
    \item \textbf{Pickup and delivery problem with time windows (PDPTW) \cite{Dumas:1991aa}:}
        The pricing problem is similar to that of VRPTW. The DIDP model uses the filter operation, a set resource variable, and the fractional knapsack expression as the dual bound function.
        We show the MP and the DIDP model for pricing based on previous work \cite{Ropke:2009aa} in Appendix~\ref{sec:additional-models}.
        We use all 56 100-customer instances from Li and Lim~\cite{Li:2001aa}.
    \item \textbf{$P||\sum w_j C_j$:}
        We generate instances following previous work~\cite{VanDenAkker:1999a}.
        In particular, we use $n = 20, 30, 40, 50$ and $m = 3, 4, 5$ with three different configurations for $p_j$ and $w_j$:
        $p_j$ uniformly sampled from $[1, 10]$ and $w_j$ uniformly sampled from $[10, 100]$,
        $p_j$ and $w_j$ uniformly sampled from $[1, 100]$,
        and $p_j$ and $w_j$ uniformly sampled from $[10, 20]$.
        For each of the 36 configurations, we generate five instances, resulting in 180 instances in total.
    \item \textbf{MRASP:}
        We use the 76 instances published in~\cite{Ghoniem:2015aa}.
\end{itemize}

We add the new features for DIDP to didp-rs v0.9.0, implemented in Rust.
%For modeling, we focus on Python, motivated by our goal of ease of modeling.
In our B\&P implementation (B\&P DIDP), we use SCIP 9.2.1 with its Python interface PySCIPOpt and DIDPPy, the Python interface for didp-rs.
For each problem class, we formulate the MIP model for the RMP using PySCIPOpt and implement a pricer by formulating the DIDP model using DIDPPy.
Using the pricer, SCIP generates variables after solving each linear relaxation in branch-and-bound.
We also implement problem-specific branching rules using PySCIPOpt (see Section~\ref{sec:case-studies} for $P||\sum w_j C_j$ and MRASP and Appendix~\ref{sec:additional-models} for VRPTW and PDPTW). 
We use Farkas pricing \cite{Gamrath:2010ab} to generate columns for infeasible linear programs, adjusting the objective value and the dual bound function in the pricing problems.
Our code is provided as supplementary material.
Each instance is solved using a single thread on an Intel Xeon Gold 6338 CPU with a time-out of 1 hour.

\subsection{Comparison with Generic Branch-and-Price Baselines}

\begin{figure}[t]
    \centering
    \scalebox{.8}{\subfile{fig_plots_benchmark.tex}}
    \caption{Percentage of instances solved over time by different approaches.\label{fig:plots_benchmark}}
\end{figure}

\begin{table}[t]
    \caption{
        Comparison of B\&P DIDP against B\&P MIP and B\&P CP in the number of pricing iterations (\#iterations), the average number of columns generated for a pricing problem (Avg. \#columns), and the average solving time for a pricing problem (Avg. time).
        We divide the competitor's value by B\&P DIDP's, averaged over the co-solved instances, with the standard deviation shown.
    }
    \centering
    \begin{tabular}{lrrr}
    \toprule
    B\&P MIP vs. B\&P DIDP & \#iterations & Avg. \#columns & Avg. time (s) \\
    \midrule
    VRPTW & $2.71\pm1.13$ & $0.26\pm0.28$ & $3.70\pm2.74$ \\
    PDPTW & $1.83\pm0.68$ & $0.56\pm0.16$ & $36.90\pm35.81$ \\
    $P||\sum w_j C_j$ & $0.56\pm0.08$ & $4.61\pm0.51$ & $199.71\pm359.81$ \\
    MRASP & $1.11\pm0.19$ & $0.66\pm0.12$ & $895.38\pm816.00$ \\
    \bottomrule
    \end{tabular}
    
    \begin{tabular}{lrrr}
    \toprule
    B\&P CP vs. B\&P DIDP & \#iterations & Avg. \#columns & Avg. time (s) \\
    \midrule
    VRPTW & $2.05\pm0.61$ & $0.25\pm0.21$ & $27.35\pm29.75$ \\
    PDPTW & $1.41\pm0.42$ & $0.56\pm0.16$ & $72.22\pm52.96$ \\
    $P||\sum w_j C_j$ & $0.50\pm0.06$ & $5.66\pm1.20$ & $174.55\pm351.63$ \\
    MRASP & $1.06\pm0.16$ & $0.89\pm0.16$ & $545.00\pm500.12$ \\
    \bottomrule
    \end{tabular}
    \label{tab:pricing-comparison}
\end{table}

We implement MIP pricing (B\&P MIP), where the pricing problems are formulated as MIP models and solved by Gurobi 13.0.2 \cite{gurobi} via its Python interface gurobipy.
We also implement CP pricing (B\&P CP), where pricing problems are formulated as CP models and solved with Google OR-Tools CP-SAT \cite{Perron:2023aa} via its Python interface.
We use an existing MIP model of pricing problems for VRPTW \cite{Irnich:2005aa} and MRASP \cite{Ghoniem:2015ab}.
We implement MIP models for the other two problem classes and CP models for all problem classes (presented in Appendix~\ref{sec:additional-models}).
We run the pricer until it finds an optimal solution or reaches the 1-hour time limit.
The DIDP, MIP, and CP solvers possibly find suboptimal solutions before the optimal one.
From these solutions, we add all variables with negative reduced cost, so different pricing solvers may generate different sets of variables and require different numbers of iterations to converge.
%This variability is known to have a significant impact on the solving trajectory.

As an automated B\&P baseline, we use GCG 3.5.5 with its Python interface PyGCGOpt.
GCG takes a compact MIP formulation as input.
We use existing MIP models for VRPTW \cite{Vigo:2014aa}, PDPTW \cite{Furtado:2017aa}, and MRASP \cite{Ghoniem:2015aa}.
For VRPTW and PDPTW, we evaluate two variations of existing MIP models \cite{Vigo:2014aa,Furtado:2017aa}, the two-index and three-index formulations, and use the better one: three-index for GCG in VRPTW and two-index for the rest.
For $P||\sum w_j C_j$, we create a compact MIP formulation, presented in Appendix~\ref{sec:additional-models}.
We also solve the compact formulations using Gurobi with gurobipy.

%\begin{table}[t]
%    \centering
%        \begin{tabular}{rrrrrrrrrr}
%        \toprule
%        Metric            & \multicolumn{3}{c}{Avg. \#iterations} & \multicolumn{3}{c}{Avg. \#columns} & \multicolumn{3}{c}{Avg. time (s)} \\
%        \cmidrule(lr){2-4} \cmidrule(lr){5-7} \cmidrule(lr){8-10}
%        Pricer            & MIP & CP & DIDP & MIP & CP & DIDP & MIP & CP & DIDP \\
%        \midrule
%        VRPTW             & $409$ & $303$ & $155$ & $2.8$ & $3.3$ & $22.7$ & $0.42$ & $2.82$ & $0.13$ \\
%        PDPTW             & $51$ & $42$ & $31$ & $4.2$ & $3.9$ & $7.3$ & $1.82$ & $2.82$ & $0.05$ \\
%        $P||\sum w_j C_j$ & $72$ & $63$ & $129$ & $4.4$ & $5.6$ & $1.0$ & $3.27$ & $4.07$ & $0.02$ \\
%        MRASP             & $49$ & $48$ & $48$ & $3.5$ & $4.9$ & $5.5$ & $22.43$ & $19.34$ & $0.23$ \\
%        \bottomrule
%    \end{tabular}
%    \caption{
%        Comparison of pricing statistics.
%        `Avg. \#iterations' is the average number of pricing iterations.
%        `Avg. \#columns' and `Avg. time' are the average number of columns and solving time in seconds for a pricing problem.
%        All metics are averaged over instances solved by all methods.
%    }
%    \label{tab:pricing-comparison}
%\end{table}

\Cref{fig:plots_benchmark} shows the percentage of optimally solved instances over time.
B\&P DIDP is stronger than the generic B\&P baselines and Gurobi on all four problems.
%At the 60-minute time limit,
%B\&P DIDP solves
%42/112 instances of the VRPTW,
%18/56 of the PDPTW,
%180/180 of the $P||\sum w_j C_j$, and
%46/76 of the MRASP.
%In contrast, B\&P MIP solves 11/112, 4/56, 65/180, and 10/76,
%while B\&P CP solves 11/112, 5/56, 57/180, and 11/76 respectively.
%GCG solves 4/112 instances of the VRPTW, 2/56 of PDPTW, and 0 for the other three problems, outperformed by Gurobi, which is inferior to B\&P DIDP.
For VRPTW, we also evaluate VRPSolver via its open-source Python interface, VRPSolverEasy~\cite{Errami:2024aa}, which
clearly shows how far behind all generic approaches are.
To our knowledge, VRPSolverEasy is not compatible with other problem classes.

Table~\ref{tab:pricing-comparison} compares B\&P DIDP with B\&P MIP and B\&P CP by showing the average relative ratio of the number of pricing iterations, the average number of generated columns for each pricing problem, and the average time to solve each pricing problem.
In VRPTW, PDPTW, and MRASP, MIP and CP require more pricing iterations than DIDP on average, possibly because fewer columns are generated per pricing iteration.
In $P||\sum w_j C_j$, DIDP always generates one column, which is an optimal solution, for each pricing problem and requires more pricing iterations than MIP and CP.
In terms of solving time, DIDP is faster than MIP and CP in all problem classes on average, up to hundreds of times, while the standard deviation is large due to the small number of co-solved instances.
Except for $P||\sum w_j C_j$, both faster solving time per iteration and the smaller number of total iterations due to more columns per iteration seem to contribute to the overall superiority of DIDP.

%Table~\ref{tab:pricing-comparison} compares B\&P MIP, B\&P CP, and B\&P DIDP in the total number of solved pricing problems, the average number of generated columns for each pricing problem, and the average time to solve each pricing problem, averaged over instances solved by all of the three methods in each problem class.
%\textbf{TODO: discuss the result in the table.}

\subsection{Ablation Study of the New DIDP Features}

\Cref{fig:plots_ablation} shows the performance of B\&P DIDP with different configurations, and Table~\ref{tab:pricing-comparison-didp} compares pricing statistics of different configurations, including the number of state expansions (the number of times line~\ref{alg:search:gen} of Algorithm~\ref{alg:labeling} is reached).
Our base configuration, `Full', uses the labeling solver with all the modeling features.
Each new feature contributes to performance improvements across multiple problem classes.
Full is best in VRPTW, but other configurations are better in other problem classes, though it does not lose much.
%Full can be used as the default option, whereas parameter tuning can improve performance.

\begin{figure}[t]
    \centering
    \scalebox{.78}{\subfile{fig_plots_ablation.tex}}
    \caption{Percentage of instances solved over time by different configurations of DIDP pricers.\label{fig:plots_ablation}}
\end{figure}

\begin{table}[t]
    \caption{
        Comparison of labeling against CAASDy and CABS in the number of pricing iterations (\#iterations), the average number of columns generated for a pricing problem (Avg. \#columns), the average solving time for a pricing problem (Avg. time), and the average number of expansions for a pricing problem (Avg. \#expansions).
        We show the competitor's value divided by Full's, averaged over the co-solved instances, with the standard deviation.
    }
    \centering
    \begin{tabular}{lrrrr}
    \toprule
    CAASDy vs. labeling & \#iterations & Avg. \#columns & Avg. time (s) & Avg. \#expansions \\
    \midrule
    VRPTW & $4.78\pm1.73$ & $0.10\pm0.13$ & $0.99\pm0.18$ & $0.87\pm0.13$ \\
    PDPTW & $3.52\pm0.79$ & $0.12\pm0.04$ & $0.63\pm0.15$ & $0.42\pm0.07$ \\
    MRASP & $2.39\pm0.80$ & $0.15\pm0.05$ & $3.92\pm8.99$ & $1.69\pm2.25$ \\
    \bottomrule
    \end{tabular}
    
    \begin{tabular}{lrrrr}
    \toprule
    CABS vs. labeling & \#iterations & Avg. \#columns & Avg. time (s) & Avg. \#expansions \\
    \midrule
    VRPTW & $1.51\pm0.50$ & $0.40\pm0.23$ & $17.93\pm21.78$ & $19.39\pm18.97$ \\
    PDPTW & $1.42\pm0.29$ & $0.44\pm0.05$ & $1.63\pm0.46$ & $1.44\pm0.43$ \\
    $P||\sum w_j C_j$ & $0.76\pm0.16$ & $2.02\pm0.34$ & $0.53\pm0.32$ & $2.66\pm0.37$ \\
    MRASP & $1.26\pm0.28$ & $0.48\pm0.09$ & $14.85\pm47.59$ & $7.89\pm13.62$ \\
    \bottomrule
    \end{tabular}
    \label{tab:pricing-comparison-didp}
\end{table}

\begin{table}[t]
    \caption{
        Comparison of Full against $-$Filter and $-$Set Res. in the number of pricing iterations (\#iterations), the average number of columns generated for a pricing problem (Avg. \#columns), the average solving time for a pricing problem (Avg. time), and the average number of expansions for a pricing problem (Avg. \#expansions).
        We show the competitor's value divided by Full's, averaged over the co-solved instances, with the standard deviation.
    }
    \centering
    \begin{tabular}{lrrrr}
    \toprule
    $-$Filter vs. Full & \#iterations & Avg. \#columns & Avg. time (s) & Avg. \#expansions \\
    \midrule
    VRPTW & $0.94\pm0.18$ & $1.07\pm0.14$ & $82.51\pm91.16$ & $27.11\pm18.23$ \\
    PDPTW & $0.94\pm0.14$ & $1.13\pm0.16$ & $4.58\pm6.53$ & $5.09\pm2.66$ \\
    MRASP & $0.93\pm0.16$ & $1.11\pm0.19$ & $29.26\pm87.70$ & $3.72\pm4.05$ \\
    \bottomrule
    \end{tabular}
    
    \begin{tabular}{lrrrr}
    \toprule
    $-$Set Res. vs. Full & \#iterations & Avg. \#columns & Avg. time (s) & Avg. \#expansions \\
    \midrule
    VRPTW & $1.01\pm0.10$ & $0.99\pm0.05$ & $11.11\pm16.36$ & $19.99\pm37.65$ \\
    PDPTW & $1.05\pm0.17$ & $0.99\pm0.07$ & $1.30\pm0.28$ & $1.42\pm0.34$ \\
    MRASP & $1.13\pm0.30$ & $0.82\pm0.10$ & $3.21\pm10.51$ & $7.29\pm21.57$ \\
    \bottomrule
    \end{tabular}
    \label{tab:pricing-comparison-didp-models}
\end{table}

\subsubsection{Comparison of Different DIDP Solvers}

We evaluate configurations using CAASDy or CABS as a DIDP solver, instead of labeling.
These algorithms expand states minimizing $f$-values.
While CAASDy is based on A*, it is generalized to support negative costs \cite{Kuroiwa:2025aa}.
In VRPTW and PDPTW, as shown in Table~\ref{tab:pricing-comparison-didp}, CAASDy reduces the number of expansions and solving time on average.
However, CAASDy also reduces the number of columns generated per iteration and results in more pricing iterations.
As a result, in Figure~\ref{fig:plots_ablation}, labeling is better in VRPTW, while CAASDy solves a few more instances in PDPTW.
In MRASP, CAASDy increases the number of expansions and solving time on average, and is outperformed by labeling.

In VRPTW, PDPTW, and MRASP, CABS increases the number of expansions and solving time compared with labeling on average.
In multiple instances of PDPTW, CABS reduces the number of expansions and solving time, resulting in a few more solved instances than labeling.
In $P||\sum w_j C_j$, the DIDP model does not have resource variables, so Full uses CAASDy, which is equivalent to labeling.
In this problem, CABS increases the number of expansions but reduces the solving time on average, which can be attributed to the low-level implementation differences between CAASDy and CABS.

In summary, expanding the state with the minimum $f$-value can solve a pricing problem faster, but overall performance may degrade due to fewer columns generated.
This result confirms the importance of heuristic pricing, which generates promising columns before exactly solving a pricing problem, as used in the literature (e.g., \cite{Ropke:2009aa,Ghoniem:2015aa}).
Designing a DIDP solver combining heuristic pricing with exact solving is future work.

\subsubsection{Comparison of Different DIDP Models }

Disabling the filter operation ($-$Filter) consistently increases the number of state expansions and the solving time, resulting in worse performance across all problem classes.
Disabling the set resource variables ($-$Set Res.) increases the number of expansions, but its effect on overall performance depends on the problem class.
In VRPTW, it significantly increases the average solving time and results in worse overall performance.
In PDPTW, the effect is moderate, but Full solves two more instances than $-$Set Res.
In MRASP, while Full is faster on average, $-$Set Res. is faster in multiple instances, resulting in a few more instances solved.
When a set resource variable is used, we may detect more dominated states, but it increases the cost of each dominance check; we must compare each state against more states sharing the same non-resource variable values, and this cost does not necessarily pay off in MRASP.

We also evaluate Profit, a configuration with a dual bound function that takes the sum over all negative reduced cost items: the fractional knapsack bound with an infinite capacity for VRPTW, PDPTW, and $P||\sum w_j C_j$ (Equation~\eqref{eq:example_bellman:bound} for VRPTW), and $\sum_{j \in \mathcal{R}} \min\{ 0, -\pi_j + w_j a_j \}$ for MRASP, ignoring the current time $t$.
The configuration No Bound does not use a dual bound function.
We show the pricing statistics in Table~\ref{tab:pricing-comparison-bound} in Appendix~\ref{sec:additional-results}.
In all problem classes, No Bound is the worst, increasing the number of expansions and solving time.
Profit also increases the number of expansions in PDPTW, $P|| \sum w_j C_j$, and MRASP.
In VRPTW and $P|| \sum w_j C_j$, Profit has almost the same solving time as Full on average.
In PDPTW, Full solves one more instance than Profit.
In MRASP, Profit increases the solving time on average, but solves the same number of instances as Full.

\section{Conclusion}
\label{sec:conclusion}

We investigate DIDP as a generic pricing solver for column generation, introducing a generic labeling algorithm and three modeling features.
Using these additions, we model the pricing problems declaratively and
implement branch-and-price for four routing and scheduling problems.
The experiments show that DIDP pricing is substantially faster than generic pricing using MIP and CP, and highlight that the new features are beneficial in multiple problem classes.
To fill the gap with problem-specific branch-and-price, developing DIDP solvers for heuristic pricing and incorporating problem-specific cuts are promising directions.

%The key takeaway is that DIDP is valuable for rapid prototyping of column generation solvers.
%The results demonstrate that
%DIDP is a practical middle ground between fully bespoke pricing codes like VRPSolver and generic constraint-based pricing,
%and
%can provide quick proof-of-concept evidence for (or against) a column generation approach
%before investing in a bespoke high-performance implementation.

%%
%% Bibliography
%%

%% Please use bibtex,

\bibliography{references}

\appendix

\section{Additional Models} \label{sec:additional-models}

We present column generation models that are not covered in the main text.
In addition, we show a compact MIP formulation for $P||\sum w_j C_j$.

\subsection{VRPTW}

In VRPTW, a directed graph $(\mathcal{N}, \mathcal{A})$ is given, where $\mathcal{N} = \{ 0, \ldots, n+1 \}$ is the set of nodes, and $\mathcal{A} \subseteq \{ \mathcal{N} \times \mathcal{N} : i \neq j, i < n + 1, j > 0 \}$.
Each customer $i$ has load $l_i > 0$, release time $a_i \ge 0$, deadline $b_i \ge 0$, and service duration $s_i \ge 0$. Each arc $(i,j)$ has distance $d_{i,j} \ge 0$ and travel cost $c_{i,j}$.
A solution is to minimize the total travel cost of a set of paths visiting all customers $\{ 1, \ldots, n \}$ exactly once.
Each path satisfies the following conditions: it starts from node $0$ and ends at node $n+1$; the cumulative load never exceeds capacity $Q$; and the visit to each node $i$ occurs within $[a_i,b_i]$.

Let $\mathcal{P}$ be a set of paths, $c_p$ be the total travel cost of path $p \in \mathcal{P}$, and $a_{p,i}$ be a constant such that $a_{p,i} = 1$ if path $p $ visits node $i \in \mathcal{N}$ and $a_{p,i} = 0$ otherwise.
\Cref{eq:cg_vrptw} shows the MP, using a binary variable $\lambda_p$ representing the selection of path $p$.
\begin{subequations}
    \label[prob]{eq:cg_vrptw}
    \begin{flalign}
        & \min \sum_{p \in \mathcal{P}} c_p \lambda_p
        \label[target]{eq:cg_vrptw:objective} \\
        & \sum_{p \in \mathcal{P}} a_{p,j} \lambda_p = 1 & \forall j = 1, \ldots, n
        \label[constr]{eq:cg_vrptw:customer-constraint} \\
        & \lambda_p \in \mathbb{Z}_+ & \forall p \in \mathcal{P}.
    \end{flalign}
\end{subequations}

For pricing, we have shown the DIDP model as a running example, using the travel cost $c_{i,j} = -\pi_j + d_{i,j}$, where $\pi_j$ is the dual value for \Cref{eq:cg_vrptw:customer-constraint}.
The MIP model is taken from previous work \cite{Irnich:2005aa}.
Branching introduces constraints specifying that edge $(j, k)$ must be or must not be used, where $(j, k)$ is the most fractional edge.

\subsubsection{CP Model for Pricing}

We use a Boolean variable $x_{i,j} \in \{ \bot, \top \}$ indicating that edge $(i, j)$ is used, a Boolean variable $y_i$ indicating that a node $i$ is visited, and an integer variable $t_i$ representing the time when node $i$ is visited.
\begin{subequations}
    \begin{align}
        \min & \sum_{i = 1, \ldots, n} -\pi_i \mathbbm{1}(y_i) + \sum_{(i,j) \in \mathcal{A}} d_{i,j} \mathbbm{1}(x_{i,j}) \\
        & \mathsf{Circuit}\left(\{ (i, j, x_{i,j}) \mid (i, j) \in \mathcal{A} \} \cup \{ (i, i, \neg y_i \mid i = 1, \ldots, n) \} \cup \{ (n+1, 0, \top) \} \right) \\
        & x_{i,j} \implies t_j \geq t_i + s_i + d_{i,j} \quad \quad \quad \quad \quad \quad \quad \quad \quad \quad \quad \quad \quad \quad \quad \quad \quad \forall (i, j) \in \mathcal{A} \\
        & \sum_{i=1, \ldots, n} l_i \mathbbm{1}(y_i) \leq Q \\
        & t_i \in [a_i, b_i] \quad \quad \quad \quad \quad \quad \quad \quad \quad \quad \quad \quad \quad \quad \quad \quad \quad \quad \quad \quad \quad \quad \quad \quad \quad \quad \forall i \in \mathcal{N} \\
        & x_{i,j} \in \{ \bot, \top \} \quad \quad \quad \quad \quad \quad \quad \quad \quad \quad \quad \quad \quad \quad \quad \quad \quad \quad \quad \quad \quad \quad \quad \forall (i, j) \in \mathcal{A} \\
        & y_i \in \{ \bot, \top \} \quad \quad \quad \quad \quad \quad \quad \quad \quad \quad \quad \quad \quad \quad \quad \quad \quad \quad \quad \quad \quad \quad \quad \forall i = 1, \ldots, n.
    \end{align}
\end{subequations}
In the objective function, we use the function $\mathbbm{1} : \{ \bot, \top \} \to \{ 0, 1 \}$ such that $\mathbbm{1}(\top) = 1$ and $\mathbbm{1}(\bot) = 0$.
We use $\mathsf{Circuit}$, the circuit global constraint \cite{Lauriere:1978aa} implemented in Google OR-Tools CP-SAT, which takes a set of triples $(i, j, x_{i,j})$ as input, ensuring that edge $(i, j)$ is used in a Hamiltonian circuit if and only if the Boolean variable $x_{i,j} = \top$.
We also introduce a triple $(i, i, \neg y_i)$, ensuring that node $i$ is not in the circuit if and only if $y_i = \bot$, and $(n + 1, 0, \top)$ to make a path from $0$ to $n + 1$ a circuit.

\subsection{PDPTW}

In PDPTW, a vehicle picks up an item at one customer and delivers it to another customer, while respecting time windows~\cite{Dumas:1991aa}.
Let $m$ be the number of tasks, $\mathcal{N} = \{1,\ldots,n\}$ be the set of tasks, and $\mathcal{L} = \{ 0, \ldots, 2n + 1 \}$  be the set of locations.
Each task $i \in \mathcal{N}$ is associated with a pickup location $i \in \mathcal{L}$ and a delivery location $n + i \in \mathcal{L}$.
Nodes $0$ and $2n+1$ represent the start and end depot locations, respectively.

Every task $i \in N$ has a load $l_i \geq 0$, and every location $i \in \mathcal{L}$ has release time $a_i \geq 0$, deadline $b_i \geq 0$ and service duration $s_i \geq 0$.
Each vehicle has the maximum capacity of $Q$.
Let $\mathcal{A} = \{ (i,j) \in \mathcal{L} \times \mathcal{L} : i \neq j, i < 2n+1 , j > 0 \}$
be the set of arcs.
Every arc $(i,j) \in A$ has a travel distance $d_{i,j} \geq 0$.
The objective is to minimize the number of vehicles used and the total travel distance.
Following the previous work~\cite{Furtado:2017aa}, our objective function is the sum of the number of used vehicles multiplied by a constant penalty $u$ and the total travel distance, where $u = 10000$.
For all evaluated methods, we tighten time windows and reduce arcs by preprocessing, following~\cite{Dumas:1991aa}.

The MP is exactly the same as \Cref{eq:cg_vrptw}.
Branching introduces constraints specifying that edge $(j, k)$ must be or must not be used.

\subsubsection{DIDP Model for Pricing}

In the pricing DIDP model, we use five state variables: $\mathcal{R}$ is the set of tasks whose pickup locations are reachable; $\mathcal{O}$ is the set of tasks whose pickup locations are visited; $i$ is the current location; $q$ is the current load; and $t$ is the current time.
The original problem corresponds to the target state $(\mathcal{R} = \mathcal{N}, \mathcal{O} = \emptyset, i = 0, q = 0, t = 0)$.
\begin{subequations}
    \label[prob]{eq:dp_pdptw}
    \begin{align}
    & \begin{array}{l}
    V(\mathcal{R}, \mathcal{O}, i, q, t) = \\
    \min\limits_{j \in \text{Next}(\mathcal{R}, \mathcal{O}, i, q, t)} \left\{ \begin{array}{ll}
        0 & \text{ if } i = 2n + 1 \\
        d_{i,j} - \pi_j + V(\mathcal{R}'(j), \mathcal{O} \cup \{ j \}, j, q + l_j, t'(j)) & \text{ if } j \in \mathcal{N} \\
        d_{i,j} + V(\mathcal{R}'(j), \mathcal{O} \setminus \{ j - n \}, j, q - l_{j-n}, t'(j)) & \text{ if } j - n \in \mathcal{N}
    \end{array} \right.
    \end{array}
    \label{eq:dp_pdptw:recursive} \\
    & \begin{array}{l}
    V(\mathcal{R}, \mathcal{O}, i, q, t) = \infty \quad \quad \quad \quad \quad \quad \quad \quad \quad \quad \quad \text{ if } \exists j \in \mathcal{O}, t + s_i + d^*_{i,n+j} > b_{n+j}
    \end{array}
    \label{eq:dp_pdptw:state-constr} \\
    & \begin{array}{l}
    V(\mathcal{R}_1, \mathcal{O}, i, q_1, t_1) \leq V(\mathcal{R}_2, \mathcal{O}, i, q_2, t_2) \quad \quad \quad \quad \text{ if } \mathcal{R}_2 \subseteq \mathcal{R}_1 \land q_1 \leq q_2 \land t_1 \leq t_2
    \end{array}
    \label[ineq]{eq:dp_pdptw:dominance} \\
    & \begin{array}{l}
    V(\mathcal{R}, \mathcal{O}, i, q, t) \geq \\
    \max\left\{
    \begin{array}{l}
        0 \text{ if } i = 2n + 1 \\
        -\textsf{fractional\_knapsack}\left(\mathcal{R}, b_{2n+1} - t - d^{\text{in}}(\mathcal{O}), \left( v^{\text{in}}_j \right)_{j = 1, \dots, n}, \left( w^{\text{in}}_j \right)_{j = 1, \dots, n} \right) \\
        -\textsf{fractional\_knapsack} \left( \mathcal{R}, b_{2n+1} - t - d^{\text{out}}(\mathcal{O}), \left( v^{\text{out}}_j \right)_{j = 1, \dots, n }, \left( w^{\text{out}}_j \right)_{j = 1, \dots, n} \right)
    \end{array} \right.
    \end{array}
    \label[ineq]{eq:dp_pdptw:bound}
    \end{align}
\end{subequations}
\Cref{eq:dp_pdptw:recursive} is the Bellman equation.
The base case is when $i = 2n+1$.
Each transition visits a location $j \in \text{Next}(\mathcal{R}, \mathcal{O}, i, q, t)$, where $\text{Next}(\mathcal{R}, \mathcal{O}, i, q, t) = \{ j \in \mathcal{R} : (i, j) \in \mathcal{A}, t + s_i + d_{i,j} \leq b_j, q + l_j \leq Q \} \cup \{ j + n : j \in \mathcal{O}, (i, j + n) \in \mathcal{A}, t + s_i + d_{i,j+n} \leq b_{j+n} \} \cup \{ 2n + 1 : q = 0 \}$.
When a location $j$ is visited, $i$ is updated to $j$, and $t$ is updated to $t'(j) = \max\{ t + s_i + d_{i,j}, a_j \}$.
Using the filter operation, $\mathcal{R}$ is updated to $\mathcal{R}'(j) = \{ k \in \mathcal{R} \setminus \{ j \} \mid t'(j) + s_j + d^*_{j,k} \leq b_k \} $, where $d^*_{j,k}$ is the precomputed shortest possible time to reach location $k$ from $j$.
If $j$ is a pickup location, $\mathcal{O}$ is updated to $\mathcal{O} \cup \{ j \}$, and $q$ is updated to $q + l_j$.
If $j$ is a delivery location, $\mathcal{O}$ is updated to $\mathcal{O} \setminus \{ j \}$, and $q$ is updated to $q - l_j$.
\Cref{eq:dp_pdptw:state-constr} ensures that the delivery location of each task in $\mathcal{O}$ must be visited by the deadline.
\Cref{eq:dp_pdptw:dominance} defines state dominance using $\mathcal{R}$, $q$, and $t$ as resource variables.
\Cref{eq:dp_pdptw:bound} is the dual bound function.
We use the fractional knapsack bound considering the deadline for the end depot.
Using $d^{\text{in}}_j = \min_{k \in \mathcal{L} : (k,j) \in \mathcal{A}} d_{k,j}$, to complete the tasks $\mathcal{O}$ and return to the end depot, we need at least $d^{\text{in}}(\mathcal{O}) = \sum_{j \in \mathcal{O}} \left( d^{\text{in}}_{n+j} + s_{n+j} \right) + d^{\text{in}}_{2n+1}$.
Completing task $j$ increases the cost by at least $v^{\text{in}}_j = d^{\text{in}}_j - \pi_j + d^{\text{in}}_{n+j}$ and the time by at least $w^{\text{in}}_j = d^{\text{in}}_j + s_j + d^{\text{in}}_{n+j} + s_{n+j}$.
We consider the 0-1 knapsack problem with the capacity $b_{2n1} - t -  d^{\text{in}}(\mathcal{O})$, where each item $j \in \mathcal{R}$ has profit $v^{\text{in}}_j$ and weight $w^{\text{in}}_j$.
We also use a similar bound using $d^{\text{out}}_j = \min_{k \in \mathcal{L} : (j,k) \in \mathcal{A}} d_{j,k}$, $d^{\text{out}}(\mathcal{O}) = d^{\text{out}}_i + \sum_{j \in \mathcal{O}} \left(s_{n+j} + d^{\text{out}}_{n+j} \right)$, $v^{\text{out}}_j = d^{\text{out}}_j - \pi_j + d^{\text{out}}_{n+j}$, and $w^{\text{out}}_j = s_j + d^{\text{out}}_j + s_{n+j} + d^{\text{out}}_{n+j}$.
The reduced cost is computed by adding the penalty $u$ to the objective value.

\subsubsection{MIP Model for Pricing}

In the pricing MIP model, we use a binary variable $x_{i,j} \in \{ 0, 1 \}$ indicating that edge $(i, j)$ is traversed, continuous variables $t_i$ and $q_i$ to represent the time and load at location $i$.
\begin{subequations}
    \begin{align}
        \min & \sum_{i \in \mathcal{N}} -\pi_i \sum_{(i, j) \in \mathcal{A}} x_{i,j} + \sum_{(i,j) \in \mathcal{A}} d_{i,j} x_{i,j} \\
        & \sum_{(0, j) \in \mathcal{A}} x_{0,j} = \sum_{(i, 2n+1) \in \mathcal{A}} x_{i,2n+1} = 1 \\
        & \sum_{(i, j) \in \mathcal{A}} x_{i,j} = \sum_{(j, i) \in \mathcal{A}} x_{j,i} & \forall i = 1, \ldots, 2n \\
        & \sum_{(i, j) \in \mathcal{A}} x_{i,j} = \sum_{(n + i, j) \in \mathcal{A}} x_{n+i,j} & \forall i \in \mathcal{N} \label[constr]{eq:mip_pdptw:pd} \\
        & t_{n+i} \geq t_i + s_i + d^*_{i,n+i}  & \forall i \in \mathcal{N} \label[constr]{eq:mip_pdptw:pd-time} \\
        & t_j \geq t_i + s_i + d_{i,j} - M (1 - x_{i,j}) & \forall (i, j) \in \mathcal{A} \label[constr]{eq:mip_pdptw:time} \\
        & q_j \geq q_i + l_j - M (1 - x_{i,j}) & \forall (i, j) \in \mathcal{A}, j \in \mathcal{N} \label[constr]{eq:mip_pdptw:load-pickup} \\
        & q_{n+j} \geq q_i - l_j - M (1 - x_{i,n+j}) & \forall (i, n+j) \in \mathcal{A}, j \in \mathcal{N} \label[constr]{eq:mip_pdptw:load-delivery} \\
        & q_0 = 0 \label[constr]{eq:mip_pdptw:q-zero} \\
        & q_i \in [l_i, Q], q_{n+i} \in [0, Q - l_i] & \forall i \in \mathcal{N} \label[constr]{eq:mip_pdptw:q} \\
        & t_i \in [a_i, b_i]  & \forall i \in \mathcal{L} \\
        & x_{i,j} \in \{ 0, 1 \} & \forall (i, j) \in \mathcal{A}.
    \end{align}
\end{subequations}
\Cref{eq:mip_pdptw:pd} ensures that the pickup location is visited if and only if the corresponding delivery location is visited.
\Cref{eq:mip_pdptw:pd-time} ensures that the delivery location is visited after the corresponding pickup location.
\Cref{eq:mip_pdptw:time} ensures that time window constraints, using $M = b_i + s_i + d_{i,j} - a_j$.
\Cref{eq:mip_pdptw:load-pickup,eq:mip_pdptw:load-delivery} ensure the capacity constraints at pickup and delivery locations, using $M = \overline{q}_i$ and $M = \overline{q}_i - l_j$, respectively, where $\overline{q}_i$ is an upper bound on $q_i$ defined in \Cref{eq:mip_pdptw:q-zero} and \eqref{eq:mip_pdptw:q}.

\subsubsection{CP Model }

The pricing CP model is based on the same idea as the MIP model and similar to that of VRPTW.
We use Boolean variable $x_{i,j} \in \{ \bot, \top \}$ to indicate that edge $(i, j)$ is traversed, $y_i \in \{ \bot, \top \}$ to indicate that task $i$ is completed, the circuit global constraint, and logical implications instead of big-M constraints.
We input triples $(i, i, \neg y_i)$ and $(n+i, n+i, \neg y_i)$ to the circuit constraint to ensure that $i$ and $n+i$ are visited if and only if $y_i = \top$.

\begin{subequations}
    \begin{align}
        \min & \sum_{i \in \mathcal{N}} -\pi_i \mathbbm{1}(y_i) + \sum_{(i,j) \in \mathcal{A}} d_{i,j} \mathbbm{1}(x_{i,j}) \\
        & \mathsf{Circuit}\left( \begin{array}{l}
            \{ (i, j, x_{i,j}) \mid (i, j) \in \mathcal{A} \} \\
            \cup \{ (i, i, \neg y_i) \mid i \in \mathcal{N}) \} \cup \{ (n+i, n+i, \neg y_i) \mid i \in \mathcal{N} \} \cup \{ (2n+1, 0, \top) \}
        \end{array} \right) \\
        & x_{i,j} \implies t_j \geq t_i + s_i + d_{i,j} \quad \quad \quad \quad \quad \quad \quad \quad \quad \quad \quad \quad \quad \quad \quad \quad \quad \forall (i, j) \in \mathcal{A} \\
        & x_{i,j} \implies q_j = q_i + l_j  \quad \quad \quad \quad \quad \quad \quad \quad \quad \quad \quad \quad \quad \quad \quad \quad \quad \forall (i, j) \in \mathcal{A}, j \in \mathcal{N} \\
        & x_{i,n+j} \implies q_{n+j} = q_i - l_j  \quad \quad \quad \quad \quad \quad \quad \quad \quad \quad \quad \quad \quad \forall (i, n+j) \in \mathcal{A}, j \in \mathcal{N} \\
        & t_{n+i} \geq t_i + s_i + d^*_{i,n+i} \quad \quad \quad \quad \quad \quad \quad \quad \quad \quad \quad \quad \quad \quad \quad \quad \quad \quad \quad \quad  \quad\forall i \in \mathcal{N} \\
        & q_0 = 0 \\
        & q_i \in [l_i, Q], q_{n+i} \in [0, Q - l_i] \quad \quad \quad \quad \quad \quad \quad \quad \quad \quad \quad \quad \quad \quad \quad \quad \quad \quad \ \forall i \in \mathcal{N} \\
        & t_i \in [a_i, b_i] \quad \quad \quad \quad \quad \quad \quad \quad \quad \quad \quad \quad \quad \quad \quad \quad \quad \quad \quad \quad \quad \quad \quad \quad \quad \quad \ \forall i \in \mathcal{L} \\
        & x_{i,j} \in \{ \bot, \top \} \quad \quad \quad \quad \quad \quad \quad \quad \quad \quad \quad \quad \quad \quad \quad \quad \quad \quad \quad \quad \quad \quad \quad \ \forall (i, j) \in \mathcal{A} \\
        & y_i \in \{ \bot, \top \} \quad \quad \quad \quad \quad \quad \quad \quad \quad \quad \quad \quad \quad \quad \quad \quad \quad \quad \quad \quad \quad \quad \quad \quad \quad \ \forall i \in \mathcal{N}.
    \end{align}
\end{subequations}

\subsection{Parallel Machine Scheduling}

We present the MIP and CP models for pricing and the compact MIP model given to GCG.

\subsubsection{MIP Model for Pricing}

In the pricing MIP model, we use a binary variable $x_j$ indicating that job $j$ is scheduled and a continuous variable $C_j$ to represent the completion time of job $j$.
\begin{subequations}
    \begin{align}
        \min & \sum_{j \in \mathcal{J}} -\pi_j x_j + w_j C_j \\
        & C_j \geq (a_j + p_j) x_j & \forall j \in \mathcal{J} \label[constr]{eq:pmwc_mip:first} \\
        & C_j \geq \sum_{k = 1}^j p_k x_k - M (1 - x_j) & \forall j \in \mathcal{J} \label[constr]{eq:pmwc_mip:time} \\
        & C_j \in [0, b_j] & \forall j \in \mathcal{J} \label[constr]{eq:pmwc_mip:c} \\
        & x_j \in \{ 0, 1 \} & \forall j \in \mathcal{J}.
    \end{align}
\end{subequations}
As mentioned in Section~\ref{sec:pmwc}, if job $j$ and job $k < j$ are scheduled, $j$ is later than $k$.
Thus, \Cref{eq:pmwc_mip:time} ensures that $C_j$ becomes at least the sum of the processing times of scheduled jobs with $k \leq j$ if $x_j = 1$, using $M = \sum_{k=1}^j p_k$.
\Cref{eq:pmwc_mip:first,eq:pmwc_mip:c} ensure the time windows.

\subsubsection{CP Model for Pricing}

The CP model is based on the same idea as the MIP model, while $x_j$ is a Boolean variable.
\begin{subequations}
    \begin{align}
        \min & \sum_{j \in \mathcal{J}} -\pi_j \mathbbm{1}(x_j) + w_j C_j \\
        & x_j \implies C_j \geq a_j + p_j \land C_j = p_j + \sum_{k = 1}^{j-1} p_k x_k  & \forall j \in \mathcal{J} \\
        & C_j \in [0, b_j] & \forall j \in \mathcal{J} \\
        & x_j \in \{ \bot, \top \} & \forall j \in \mathcal{J}.
    \end{align}
\end{subequations}

\subsubsection{Compact MIP Model}

In the compact MIP model, binary variable $x_{j,i}$ indicates that job $j$ is assigned to machine $i$.
\begin{subequations}
    \begin{align}
        \min & \sum_{j \in \mathcal{J}} w_j C_j \\
        & \sum_{i = 1}^m x_{j,i} = 1 & \forall j \in \mathcal{J} \label[constr]{eq:pmwc_compact:machine} \\
        & C_j \geq \sum_{k=1}^j p_k x_{k,i} - M (1 - x_{j,i}) & \forall j \in \mathcal{J}, \forall i=1,\ldots,m  \\
        & C_j \in [a_j + p_j, b_j] & \forall j \in \mathcal{J} \\
        & x_{j,i} \in \{ 0, 1 \} & \forall j \in \mathcal{J}, \forall i = 1, \ldots, m.
    \end{align}
\end{subequations}
\Cref{eq:pmwc_compact:machine} ensures a job is processed by exactly one machine.

\subsection{CP Model for MRASP}

The pricing CP model is similar to the existing pricing MIP model \cite{Ghoniem:2015aa} but uses logical constraints.
Boolean variable $x_i$ indicates that aircraft $i$ is scheduled, Boolean variable $y_{i,j}$ indicates that aircraft $j$ is scheduled later than $i$, and an integer variable $t_i$ represents the scheduled time for aircraft $i$.
\begin{subequations}
    \begin{align}
        \min & \sum_{i \in \mathcal{N}} -\pi_i \mathbbm{1}(x_i) +  w_i t_i  \\
        & x_i \implies t_i \geq a_i & \forall i \in \mathcal{N} \label[constr]{eq:mrasp_cp:release} \\
        & y_{i,j} \implies x_i \land t_j \geq t_i + d_{g_i,o_i,g_j,o_j} & \forall i \in \mathcal{N}, \forall j \in \mathcal{N} \setminus \{ i \} \label[constr]{eq:mrasp_cp:separation} \\
        & \textsf{AtMostOne}(\{ y_{i,j}, y_{j,i} \}) & \forall i \in \mathcal{N}, \forall j \in \mathcal{N} \setminus \{ i \} \label[constr]{eq:mrasp_cp:at_most_one} \\
        & x_i \land x_j \implies y_{i,j} \lor y_{j,i} & \forall i \in \mathcal{N}, \forall j \in \mathcal{N} \setminus \{ i \} \label[constr]{eq:mrasp_cp:or} \\
        & \neg y_{i,j}  & \forall i \in \mathcal{N}, \forall j \in \mathcal{N} \setminus \{ i \}, a_i + d_{g_i,o_i,g_j,o_j} > b_j \label[constr]{eq:mrasp_cp:fix1} \\
        & \neg y_{j,i} & \begin{array}{l}
            \forall i \in \mathcal{N}, \forall j \in \mathcal{N} \setminus \{ i \} \\
            a_i < a_j, b_i \leq b_j, g_i = g_j, o_i = o_j
        \end{array} \label[constr]{eq:mrasp_cp:fix2} \\
        & t_i \in [0, b_i] & \forall i \in \mathcal{N} \label[constr]{eq:mrasp_cp:t} \\
        & x_i \in \{ \bot, \top \} & \forall i \in \mathcal{N} \\
        & y_{i,j} \in \{ \bot, \top \} & \forall i \in \mathcal{N}, \forall j \in \mathcal{N} \setminus \{ i \}.
    \end{align}
\end{subequations}
\Cref{eq:mrasp_cp:release,eq:mrasp_cp:t} ensure the time window constraints.
\Cref{eq:mrasp_cp:separation} ensures the separation time and link $x_i$ and $y_{i,j}$ together with \Cref{eq:mrasp_cp:or}  .
\Cref{eq:mrasp_cp:fix1,eq:mrasp_cp:fix2} exclude impossible or suboptimal precedences and are adapted from the pricing MIP model used in the previous work.

In addition, in the root node of B\&P, we introduce the following constraint, excluding aircraft $i$ that cannot be scheduled with a negative reduced cost:
\begin{subequations}
    \begin{align}
        & \neg x_i & \forall i \in \mathcal{N}, -\pi_i + w_i a_i \geq 0.
    \end{align}
\end{subequations}
An equivalent constraint is also introduced in the pricing MIP model.

In branching, constraints to forbid or enforce that aircraft $j$ is scheduled directly after $i$ are introduced.
Let $\mathcal{A}_\text{forbidden}$ and $\mathcal{A}_\text{enforced}$ be the set of pairs $(i, j)$ such that scheduling $j$ directly after $i$ is forbidden or enforced.
Then, we introduce the following constraints:
\begin{subequations} \label[constr]{eq:mrasp_cp:branch}
    \begin{align}
        & x_i \implies \bigvee_{j \in \mathcal{N} \setminus \{ i \}} y_{j,i} & \forall (0, i) \in \mathcal{A}_\text{forbidden} \\
        & x_i \implies \bigvee_{j \in \mathcal{N} \setminus \{ i \}} y_{i,j} & \forall (i, n+1) \in \mathcal{A}_\text{forbidden} \\
        & y_{i,j} \implies \sum_{k \in \mathcal{N} \setminus \{ j \}} y_{k,j} - \sum_{k \in \mathcal{N} \setminus \{ i \}} y_{k,i} \geq 2   & \forall (i, j)\in \mathcal{A}_\text{forbidden} \label[constr]{eq:mrasp_cp:forbid} \\
        & \neg y_{j,i} & \forall (0, i) \in \mathcal{A}_{\text{enforced}}, \forall j \in \mathcal{N} \setminus \{ i \} \\
        & \neg y_{i,j} & \forall (i, n+1) \in \mathcal{A}_{\text{enforced}}, \forall j \in \mathcal{N} \setminus \{ i \} \\
        & x_i \implies y_{i,j} & \forall (i, j) \in \mathcal{A}_{\text{enforced}} \\
        & x_j \implies y_{i,j} & \forall (i, j) \in \mathcal{A}_{\text{enforced}} \\
        & y_{i,j} \implies \sum_{k \in \mathcal{N} \setminus \{ j \}} y_{k,j} - \sum_{k \in \mathcal{N} \setminus \{ i \}} y_{k,i} = 1  & \forall (i, j) \in \mathcal{A}_{\text{enforced}} \label[constr]{eq:mrasp_cp:enforce}.
    \end{align}
\end{subequations}
\Cref{eq:mrasp_cp:forbid} ensures that at least one aircraft is scheduled after $i$ and before $j$.
In a similar way, \Cref{eq:mrasp_cp:enforce} ensures that no aircraft is scheduled after $i$ and before $j$.
The linearized versions of the above constraints are used in the pricing MIP model.

\section{Additional Results} \label{sec:additional-results}

We present the pricing statistics for B\&P DIDP configurations with different dual bound functions in Table~\ref{tab:pricing-comparison-bound}.

\begin{table}[h]
    \caption{
        Comparison of Full against Profit and No Bound in the number of pricing iterations (\#iterations), the average number of columns generated for a pricing problem (Avg. \#columns), the average solving time for a pricing problem (Avg. time), and the average number of expansions for a pricing problem (Avg. \#expansions).
        We show the competitor's value divided by Full's, averaged over the co-solved instances, with the standard deviation.
    }
    \centering
    \begin{tabular}{lrrrr}
    \toprule
    Profit vs. Full & \#iterations & Avg. \#columns & Avg. time (s) & Avg. \#expansions \\
    \midrule
    VRPTW & $1.00\pm0.00$ & $1.00\pm0.00$ & $0.97\pm0.04$ & $1.00\pm0.00$ \\
    PDPTW & $1.00\pm0.00$ & $1.00\pm0.00$ & $1.24\pm0.45$ & $1.35\pm0.52$ \\
    $P||\sum w_j C_j$ & $1.01\pm0.17$ & $1.00\pm0.01$ & $1.00\pm0.06$ & $1.10\pm0.06$ \\
    MRASP & $0.98\pm0.18$ & $1.01\pm0.07$ & $2.72\pm6.10$ & $1.68\pm1.90$ \\
    \bottomrule
    \end{tabular}
    
    % Table 3: Dual bound configuration No Bound
    \begin{tabular}{lrrrr}
    \toprule
    No Bound vs. Full & \#iterations & Avg. \#columns & Avg. time (s) & Avg. \#expansions \\
    \midrule
    VRPTW & $1.00\pm0.00$ & $1.00\pm0.00$ & $1.12\pm0.39$ & $1.11\pm0.11$ \\
    PDPTW & $1.00\pm0.00$ & $1.00\pm0.00$ & $1.60\pm0.74$ & $2.16\pm1.04$ \\
    $P||\sum w_j C_j$ & $0.31\pm0.10$ & $21.14\pm11.86$ & $12.18\pm16.57$ & $94.71\pm119.91$ \\
    MRASP & $1.08\pm0.42$ & $0.98\pm0.08$ & $11.19\pm24.76$ & $2.37\pm2.71$ \\
    \bottomrule
    \end{tabular}
        \label{tab:pricing-comparison-bound}
    \end{table}
\end{document}

%% file: fig_plots_benchmark.tex
\begin{tikzpicture}

    \pgfplotstableread[col sep=comma]{results/vrptw_labeling-knapsack.csv}{\vrptwlabeling}
    \pgfplotstableread[col sep=comma]{results/vrptw_mip.csv}{\vrptwmip}
    \pgfplotstableread[col sep=comma]{results/vrptw_cp.csv}{\vrptwcp}
    \pgfplotstableread[col sep=comma]{results/vrptw_gcg-threeindex.csv}{\vrptwgcgthreeindex}
    %\pgfplotstableread[col sep=comma]{results/vrptw_gcg-twoindex.csv}{\vrptwgcgtwoindex}
    \pgfplotstableread[col sep=comma]{results/vrptw_gurobi-twoindex.csv}{\vrptwgurobitwoindex}
    % \pgfplotstableread[col sep=comma]{results/vrptw_scip-twoindex.csv}{\vrptwsciptwoindex}
    \pgfplotstableread[col sep=comma]{results/vrptw_vrpse.csv}{\vrptwvrpse}

    \pgfplotstableread[col sep=comma]{results/pdptw_labeling-knapsack.csv}{\pdptwlabeling}
    \pgfplotstableread[col sep=comma]{results/pdptw_mip.csv}{\pdptwmip}
    \pgfplotstableread[col sep=comma]{results/pdptw_cp.csv}{\pdptwcp}
    %\pgfplotstableread[col sep=comma]{results/pdptw_gcg-threeindex.csv}{\pdptwgcgthreeindex}
    \pgfplotstableread[col sep=comma]{results/pdptw_gcg-twoindex.csv}{\pdptwgcgtwoindex}
    \pgfplotstableread[col sep=comma]{results/pdptw_gurobi-twoindex.csv}{\pdptwgurobitwoindex}
    % \pgfplotstableread[col sep=comma]{results/pdptw_scip-twoindex.csv}{\pdptwsciptwoindex}

    \pgfplotstableread[col sep=comma]{results/pmwc_labeling-knapsack.csv}{\pmwclabeling}
    \pgfplotstableread[col sep=comma]{results/pmwc_mip.csv}{\pmwcmip}
    \pgfplotstableread[col sep=comma]{results/pmwc_cp.csv}{\pmwccp}
    \pgfplotstableread[col sep=comma]{results/pmwc_gcg.csv}{\pmwcgcg}
    \pgfplotstableread[col sep=comma]{results/pmwc_gurobi.csv}{\pmwcgurobi}
    % \pgfplotstableread[col sep=comma]{results/pmwc_scip.csv}{\pmwcscip}

    \pgfplotstableread[col sep=comma]{results/mrasp_labeling-reduce.csv}{\mrasplabeling}
    \pgfplotstableread[col sep=comma]{results/mrasp_mip.csv}{\mraspmip}
    \pgfplotstableread[col sep=comma]{results/mrasp_cp.csv}{\mraspcp}
    \pgfplotstableread[col sep=comma]{results/mrasp_gcg.csv}{\mraspgcg}
    \pgfplotstableread[col sep=comma]{results/mrasp_gurobi.csv}{\mraspgurobi}
    % \pgfplotstableread[col sep=comma]{results/mrasp_scip.csv}{\mraspscip}

    \begin{groupplot}[
        group style = {
            %group size=2 by 2,
            group size=4 by 1,
            %horizontal sep=13ex,
            horizontal sep=1ex,
            %vertical sep=11ex,
            xlabels at=edge bottom,
            ylabels at=edge left,
            xticklabels at=edge bottom,
            yticklabels at=edge left,
        },
        success_rate_subplot,
        ]

        \nextgroupplot[
            title={VRPTW},
            legend style={legend columns=6},
            legend to name=legend1,
            % reverse legend,
        ]
        \addplot[plotline solid, const plot mark left, mark size=0.15ex, black]  table [x=Time, y=Percentage] {\vrptwlabeling};
        \addplot[plotline dashed, const plot mark left, mark size=0.15ex, red]    table [x=Time, y=Percentage] {\vrptwmip};
        \addplot[plotline dotted, const plot mark left, mark size=0.15ex, green]  table [x=Time, y=Percentage] {\vrptwcp};
        \addplot[plotline dashdot, const plot mark left, mark size=0.15ex, blue]   table [x=Time, y=Percentage] {\vrptwgcgthreeindex};
        %\addplot[plotline loosely dotted, const plot mark left, mark size=0.15ex, purple] table [x=Time, y=Percentage] {\vrptwgcgtwoindex};
        \addplot[plotline longdash, const plot mark left, mark size=0.15ex, orange] table [x=Time, y=Percentage] {\vrptwgurobitwoindex};
        \addplot[plotline dashdotdot, const plot mark left, mark size=0.15ex, pink]   table [x=Time, y=Percentage] {\vrptwvrpse};
        \addlegendentry{B\&P DIDP};
        \addlegendentry{B\&P MIP};
        \addlegendentry{B\&P CP};
        \addlegendentry{GCG};
        %\addlegendentry{GCG -- Three-Index};
        %\addlegendentry{GCG -- Two-Index};
        \addlegendentry{Gurobi};
        % \addlegendentry{SCIP};
        \addlegendentry{VRPSolverEasy};

        \nextgroupplot[
            title={PDPTW},
            %ymax=0.4,
            %ytick distance=0.1,
        ]
        \addplot[plotline solid, const plot mark left, mark size=0.15ex, black]  table [x=Time, y=Percentage] {\pdptwlabeling};
        \addplot[plotline dashed, const plot mark left, mark size=0.15ex, red]    table [x=Time, y=Percentage] {\pdptwmip};
        \addplot[plotline dotted, const plot mark left, mark size=0.15ex, green]  table [x=Time, y=Percentage] {\pdptwcp};
        %\addplot[plotline dashdot, const plot mark left, mark size=0.15ex, blue]   table [x=Time, y=Percentage] {\pdptwgcgthreeindex};
        \addplot[plotline dashdot, const plot mark left, mark size=0.15ex, blue] table [x=Time, y=Percentage] {\pdptwgcgtwoindex};
        \addplot[plotline longdash, const plot mark left, mark size=0.15ex, orange] table [x=Time, y=Percentage] {\pdptwgurobitwoindex};

        \nextgroupplot[
            title={$P||\sum w_j C_j$},
        ]
        \addplot[plotline solid, const plot mark left, mark size=0.15ex, black]  table [x=Time, y=Percentage] {\pmwclabeling};
        \addplot[plotline dashed, const plot mark left, mark size=0.15ex, red]    table [x=Time, y=Percentage] {\pmwcmip};
        \addplot[plotline dotted, const plot mark left, mark size=0.15ex, green]  table [x=Time, y=Percentage] {\pmwccp};
        \addplot[plotline dashdot, const plot mark left, mark size=0.15ex, blue]   table [x=Time, y=Percentage] {\pmwcgcg};
        \addplot[plotline longdash, const plot mark left, mark size=0.15ex, orange] table [x=Time, y=Percentage] {\pmwcgurobi};

        \nextgroupplot[
            title={MRASP},
        ]
        \addplot[plotline solid, const plot mark left, mark size=0.15ex, black]  table [x=Time, y=Percentage] {\mrasplabeling};
        \addplot[plotline dashed, const plot mark left, mark size=0.15ex, red]    table [x=Time, y=Percentage] {\mraspmip};
        \addplot[plotline dotted, const plot mark left, mark size=0.15ex, green]  table [x=Time, y=Percentage] {\mraspcp};
        \addplot[plotline dashdot, const plot mark left, mark size=0.15ex, blue]   table [x=Time, y=Percentage] {\mraspgcg};
        \addplot[plotline longdash, const plot mark left, mark size=0.15ex, orange] table [x=Time, y=Percentage] {\mraspgurobi};

    \end{groupplot}

    % Legend
    %\coordinate (legend1_left) at (group c1r2.south west);
    %\coordinate (legend1_right) at ($(group c2r2.south east |- group c1r2.south west)$);
    %\coordinate (legend1_center) at ($(legend1_left)!0.5!(legend1_right)$);
    %\node[below=7ex] at (legend1_center) {\pgfplotslegendfromname{legend1}};
    \coordinate (legend1_center) at ($(group c1r1.south west)!0.5!(group c4r1.south east)$);
    \node[below=7ex] at (legend1_center) {\pgfplotslegendfromname{legend1}};

\end{tikzpicture}

%% file: fig_plots_ablation.tex
\begin{tikzpicture}

    \pgfplotstableread[col sep=comma]{results/vrptw_labeling-knapsack-disable-set-resource.csv}{\vrptwlabelingknapsackdisablesetresource}
    \pgfplotstableread[col sep=comma]{results/vrptw_labeling-knapsack-disable-filter.csv}{\vrptwlabelingknapsackdisablefilter}
    \pgfplotstableread[col sep=comma]{results/vrptw_caasdy-knapsack.csv}{\vrptwcaasdyknapsack}
    \pgfplotstableread[col sep=comma]{results/vrptw_cabs-knapsack.csv}{\vrptwcabsknapsack}
    \pgfplotstableread[col sep=comma]{results/vrptw_labeling-none.csv}{\vrptwlabelingnone}
    \pgfplotstableread[col sep=comma]{results/vrptw_labeling-profit.csv}{\vrptwlabelingprofit}
    \pgfplotstableread[col sep=comma]{results/vrptw_labeling-knapsack.csv}{\vrptwlabelingknapsack}

    \pgfplotstableread[col sep=comma]{results/pdptw_labeling-knapsack-disable-set-resource.csv}{\pdptwlabelingknapsackdisablesetresource}
    \pgfplotstableread[col sep=comma]{results/pdptw_labeling-knapsack-disable-filter.csv}{\pdptwlabelingknapsackdisablefilter}
    \pgfplotstableread[col sep=comma]{results/pdptw_caasdy-knapsack.csv}{\pdptwcaasdyknapsack}
    \pgfplotstableread[col sep=comma]{results/pdptw_cabs-knapsack.csv}{\pdptwcabsknapsack}
    \pgfplotstableread[col sep=comma]{results/pdptw_labeling-none.csv}{\pdptwlabelingnone}
    \pgfplotstableread[col sep=comma]{results/pdptw_labeling-profit.csv}{\pdptwlabelingprofit}
    \pgfplotstableread[col sep=comma]{results/pdptw_labeling-knapsack.csv}{\pdptwlabelingknapsack}

    \pgfplotstableread[col sep=comma]{results/pmwc_caasdy-knapsack.csv}{\pmwccaasdyknapsack}
    \pgfplotstableread[col sep=comma]{results/pmwc_cabs-knapsack.csv}{\pmwccabsknapsack}
    \pgfplotstableread[col sep=comma]{results/pmwc_labeling-none.csv}{\pmwclabelingnone}
    \pgfplotstableread[col sep=comma]{results/pmwc_labeling-profit.csv}{\pmwclabelingprofit}
    \pgfplotstableread[col sep=comma]{results/pmwc_labeling-knapsack.csv}{\pmwclabelingknapsack}

    \pgfplotstableread[col sep=comma]{results/mrasp_labeling-reduce-disable-set-resource.csv}{\mrasplabelingreducedisablesetresource}
    \pgfplotstableread[col sep=comma]{results/mrasp_labeling-reduce-disable-filter.csv}{\mrasplabelingreducedisablefilter}
    \pgfplotstableread[col sep=comma]{results/mrasp_caasdy-reduce.csv}{\mraspcaasdyreduce}
    \pgfplotstableread[col sep=comma]{results/mrasp_cabs-reduce.csv}{\mraspcabsreduce}
    \pgfplotstableread[col sep=comma]{results/mrasp_labeling-none.csv}{\mrasplabelingnone}
    \pgfplotstableread[col sep=comma]{results/mrasp_labeling-profit.csv}{\mrasplabelingprofit}
    \pgfplotstableread[col sep=comma]{results/mrasp_labeling-reduce.csv}{\mrasplabelingknapsack}

    \begin{groupplot}[
        group style = {
            group size=4 by 1,
            horizontal sep=1ex,
            %horizontal sep=13ex,
            %vertical sep=11ex,
            % xlabels at=edge bottom,
            % ylabels at=edge left,
            % xticklabels at=edge bottom,
            % yticklabels at=edge left,
        },
        success_rate_subplot,
        ]

        \nextgroupplot[
            title={VRPTW},
            ymax=0.4,
            ytick distance=0.1,
            legend style={legend columns=7},
            legend to name=legend2,
            reverse legend,
        ]
        \addplot[plotline dotted, const plot mark left, mark size=0.15ex, green]  table [x=Time, y=Percentage] {\vrptwlabelingnone};
        \addplot[plotline dashed, const plot mark left, mark size=0.15ex, red]    table [x=Time, y=Percentage] {\vrptwlabelingprofit};
        \addplot[plotline dashdotdot, const plot mark left, mark size=0.15ex, pink] table [x=Time, y=Percentage] {\vrptwlabelingknapsackdisablesetresource};
        \addplot[plotline loosely dotted, const plot mark left, mark size=0.15ex, purple]   table [x=Time, y=Percentage] {\vrptwlabelingknapsackdisablefilter};
        \addplot[plotline dashdot, const plot mark left, mark size=0.15ex, blue]   table [x=Time, y=Percentage] {\vrptwcabsknapsack};
        \addplot[plotline longdash, const plot mark left, mark size=0.15ex, orange] table [x=Time, y=Percentage] {\vrptwcaasdyknapsack};
        \addplot[plotline solid, const plot mark left, mark size=0.15ex, black]  table [x=Time, y=Percentage] {\vrptwlabelingknapsack};
        \addlegendentry{No Bound};
        \addlegendentry{Profit};
        \addlegendentry{$-$Set Res.};
        \addlegendentry{$-$Filter};
        \addlegendentry{CABS};
        \addlegendentry{CAASDy};
        \addlegendentry{Full};

        \nextgroupplot[
            title={PDPTW},
            ymax=0.4,
            ytick distance=0.1,
            yticklabels=\empty,
            ylabel={}
        ]
        \addplot[plotline dotted, const plot mark left, mark size=0.15ex, green]  table [x=Time, y=Percentage] {\pdptwlabelingnone};
        \addplot[plotline dashed, const plot mark left, mark size=0.15ex, red]    table [x=Time, y=Percentage] {\pdptwlabelingprofit};
        \addplot[plotline dashdotdot, const plot mark left, mark size=0.15ex, pink] table [x=Time, y=Percentage] {\pdptwlabelingknapsackdisablesetresource};
        \addplot[plotline loosely dotted, const plot mark left, mark size=0.15ex, purple]   table [x=Time, y=Percentage] {\pdptwlabelingknapsackdisablefilter};
        \addplot[plotline dashdot, const plot mark left, mark size=0.15ex, blue]   table [x=Time, y=Percentage] {\pdptwcabsknapsack};
        \addplot[plotline longdash, const plot mark left, mark size=0.15ex, orange] table [x=Time, y=Percentage] {\pdptwcaasdyknapsack};
        \addplot[plotline solid, const plot mark left, mark size=0.15ex, black]  table [x=Time, y=Percentage] {\pdptwlabelingknapsack};

        \nextgroupplot[
            title={$P||\sum w_j C_j$},
            ymin=0.4,
            ymax=1.0,
            ytick distance=0.2,
            xmin=0, xmax=20, xtick={0, 5, 10, 15, 20},
            ylabel={}
        ]
        \addplot[plotline dotted, const plot mark left, mark size=0.15ex, green]  table [x=Time, y=Percentage] {\pmwclabelingnone};
        \addplot[plotline dashed, const plot mark left, mark size=0.15ex, red]    table [x=Time, y=Percentage] {\pmwclabelingprofit};
        \addplot[plotline dashdot, const plot mark left, mark size=0.15ex, blue]   table [x=Time, y=Percentage] {\pmwccabsknapsack};
        \addplot[plotline longdash, const plot mark left, mark size=0.15ex, orange] table [x=Time, y=Percentage] {\pmwccaasdyknapsack};
        \addplot[plotline solid, const plot mark left, mark size=0.15ex, black]  table [x=Time, y=Percentage] {\pmwclabelingknapsack};

        \pgfplotsset{every axis/.append style={xshift=4.5ex}}

        \nextgroupplot[
            title={MRASP},
            ymin=0.4,
            ymax=1.0,
            ytick distance=0.2,
            yticklabels=\empty,
            ylabel={}
        ]
        \addplot[plotline dotted, const plot mark left, mark size=0.15ex, green]  table [x=Time, y=Percentage] {\mrasplabelingnone};
        \addplot[plotline dashed, const plot mark left, mark size=0.15ex, red]    table [x=Time, y=Percentage] {\mrasplabelingprofit};
        \addplot[plotline dashdotdot, const plot mark left, mark size=0.15ex, pink] table [x=Time, y=Percentage] {\mrasplabelingreducedisablesetresource};
        \addplot[plotline loosely dotted, const plot mark left, mark size=0.15ex, purple]   table [x=Time, y=Percentage] {\mrasplabelingreducedisablefilter};
        \addplot[plotline dashdot, const plot mark left, mark size=0.15ex, blue]   table [x=Time, y=Percentage] {\mraspcabsreduce};
        \addplot[plotline longdash, const plot mark left, mark size=0.15ex, orange] table [x=Time, y=Percentage] {\mraspcaasdyreduce};
        \addplot[plotline solid, const plot mark left, mark size=0.15ex, black]  table [x=Time, y=Percentage] {\mrasplabelingknapsack};

    \end{groupplot}

    % Legend
    %\coordinate (legend1_left) at (group c1r2.south west);
    %\coordinate (legend1_right) at ($(group c2r2.south east |- group c1r2.south west)$);
    %\coordinate (legend1_center) at ($(legend1_left)!0.5!(legend1_right)$);
    %\node[below=7ex] at (legend1_center) {\pgfplotslegendfromname{legend1}};
    \coordinate (legend2_center) at ($(group c1r1.south west)!0.5!(group c4r1.south east)$);
    \node[below=7ex] at (legend2_center) {\pgfplotslegendfromname{legend2}};

\end{tikzpicture}